\documentclass{amsart}
\usepackage[francais]{babel}
\usepackage{amsmath}
\usepackage{amssymb}

\newtheorem{thm}{Th\'eor\`eme}[section]
\newtheorem{prop}[thm]{Proposition}
\newtheorem*{prop*}{Proposition}
\newtheorem{lem}[thm]{Lemme}
\newtheorem{cor}[thm]{Corollaire}
\newtheorem*{cor*}{Corollaire}
\theoremstyle{definition}
\newtheorem{defn}[thm]{D\'efinition}
\theoremstyle{remark}
\newtheorem{rem}[thm]{Remarque}
\newtheorem{exm}[thm]{Exemple}

\newcommand{\cD}{\mathcal{D}}
\newcommand{\cH}{\mathcal{H}}
\newcommand{\cO}{\mathcal{O}}
\newcommand{\cP}{\mathcal{P}}
\newcommand{\cS}{\mathcal{S}}
\newcommand{\fS}{\mathfrak{S}}
\newcommand{\gt}{\mathfrak{t}}

\newcommand{\bb}{\mathbf{b}}
\newcommand{\bd}{\mathbf{d}}
\newcommand{\bm}{\mathbf{m}}
\newcommand{\bn}{\mathbf{n}}
\newcommand{\bs}{\mathbf{s}}
\newcommand{\bu}{\mathbf{u}}

\newcommand{\tWL}{\tilde W_L}
\newcommand{\balpha}{{\boldsymbol{\alpha}}}
\newcommand{\bbeta}{{\boldsymbol{\beta}}}
\newcommand{\bgamma}{{\boldsymbol{\gamma}}}
\newcommand{\bmu}{{\boldsymbol{\mu}}}
\newcommand{\tbalpha}{\tilde\balpha}

\newcommand{\refl}{\text{{\rm r\'efl}}}

\newcommand{\plein}{\text{{\rm plein}}}
\newcommand{\spets}{\text{{\rm spets}}}
\newcommand{\GL}{\operatorname{GL}}

\newcommand{\Ind}{\operatorname{Ind}}
\newcommand{\Irr}{\operatorname{Irr}}

\newcommand{\Gal}{\operatorname{Gal}}
\newcommand{\rot}{\operatorname{r}}
\newcommand{\Tr}{\operatorname{Tr}}
\newcommand{\bigboxtimes}{\mathop{\vcenter{\hbox{\fontsize{17.28}{0pt}%
  \selectfont$\boxtimes$}}}} 

\newcommand{\C}{\mathbb{C}}
\newcommand{\Q}{\mathbb{Q}}
\newcommand{\N}{\mathbb{N}}
\newcommand{\Z}{\mathbb{Z}}

\newcommand{\ie}{{\it i.e.,\ }}

\title[Repr\'esentations de Springer]{Repr\'esentations de Springer pour les groupes de r\'eflexions complexes imprimitifs}
\author{Pramod N.~Achar et Anne-Marie Aubert}
\thanks{Le premier auteur \'etait partiellement appuy\'e par la subvention 
DMS-0500873 de la NSF}
\address{Department of Mathematics\\
  Louisiana State University\\
  Baton Rouge, LA 70803, USA}
\email{pramod@math.lsu.edu}
\address{Institut de Math\'ematiques de Jussieu\\
UMR 7586 du C.N.R.S.\\
F-75252 Paris Cedex 05
  \\
  France}
\email{aubert@math.jussieu.fr}

\begin{document}

\begin{abstract}
\`A un groupe de r\'eflexions complexe sp\'etsial, muni d'un r\'eseau radiciel au sens de Nebe, nous associons un certain ensemble fini qui doit jouer un r\^ole analogue \`a celui de l'ensemble des classes unipotentes d'un groupe alg\'ebrique. Dans le cas des groupes imprimitifs, nous en donnons un param\'etrage combinatoire en termes des symboles g\'en\'eralis\'es de Malle et Shoji. Ce r\'esultat fournit un lien entre les travaux de Shoji sur les fonctions de Green pour les groupes de r\'eflexions complexes et ceux de Brou\'e, Kim, Malle, Rouquier, {\it et al.} sur les alg\`ebres de Hecke cyclotomiques et leurs familles de caract\`eres.
\end{abstract}

\maketitle

%**************************************************************************
\section{Introduction}
%**************************************************************************

Les groupes de r\'eflexions complexes, et surtout ceux dits sp\'etsiaux,
se sont r\'ecemment montr\'es proches des groupes de Weyl des groupes
alg\'ebriques dans de nombreux aspects: ils admettent des alg\`ebres de
Hecke et des groupes de tresses avec de bonnes propri\'et\'es; leurs
caract\`eres se r\'epartissent en ``familles''; et pour certains d'entre
eux --- les groupes imprimitifs --- Shoji a d\'evelopp\'e une th\'eorie de
fonctions de Green \cite{Sh1, Sh2, Sh3}.

Rappelons que dans le cadre des groupes alg\'ebriques r\'eductifs sur un
corps fini, les fonctions de Green sont certaines fonctions \`a valeurs
complexes d\'efinies sur l'ensemble des \'el\'ements unipotents. Elles se
calculent par un algorithme, d\^u \`a Lusztig et Shoji, qui ne d\'epend
que du groupe de Weyl. Une question naturelle est donc: est-il possible 
d'effectuer le m\^eme algorithme pour les groupes de r\'eflexions complexes?
Cette question est le point de d\'epart des travaux de Shoji, et il a
d\'ecouvert que les nouvelles ``fonctions de Green'' ainsi obtenues
semblent v\'erifier certaines conditions remarquables d'int\'egralit\'e et
de positivit\'e (en commun avec les ``vraies'' fonctions de Green),
bien qu'elles n'aient pas (encore?) d'interpr\'etation g\'eom\'etrique.

Cependant, pour d\'emarrer l'algorithme pour les groupes de Weyl, il faut
d'abord connaitre la correspondance de Springer. Par contre, pour les
groupes de r\'eflexions complexes, puisqu'il n'y a  ni vari\'et\'e
unipotente, ni correspondance de Springer, il faut choisir et imposer
sur l'ensemble de repr\'esentations irr\'eductibles une structure qui
ressemble \`a celles provenant des correspondances de Springer. Shoji n'a
trait\'e que les groupes de r\'eflexions imprimitifs, et il a choisi une
structure d\'efinie en termes des objets combinatoires dits ``symboles.'' 

Cette circonstance donne lieu a plusieurs questions: Est-ce que la
structure choisie par Shoji est pr\'ef\'er\'ee ou naturelle en un certain
sens, ou bien,  est-ce que d'autres choix donneraient lieu \`a des
fonctions de Green diff\'erentes de celles de Shoji mais \'egalement
valables? D'autre part, comment peut-on \'etendre ses r\'esultats au cas
primitif (\ie exceptionnel), o\`u l'on ne peut pas utiliser d'objets
combinatoires?

Le but de cet article est d'essayer de r\'epondre \`a ces questions. Nous
proposons ici une nouvelle construction alg\'ebrique qui associe \`a un
chaque groupe de r\'eflexions complexes (muni d'un r\'eseau radiciel) un
ensemble qui doit jouer le r\^ole de l'ensemble de classes unipotentes. Le
r\'esultat principal affirme que les symboles de Shoji sont compatibles
dans un certain sens avec notre construction, et donc que son choix
\'etait bien naturel. D'autre part, notre contruction fonctionne
\'egalement bien pour tous les groupes de r\'eflexions complexes
sp\'etsiaux, et nous obtenons ainsi les d\'ebuts d'une extension des
travaux de Shoji aux groupes primitifs.

Nous commen\c cons \`a la Section~\ref{sect:symboles} par d\'efinir tous
les objets combinatoires dont nous aurons besoin. La
section~\ref{sect:groupe} est consacr\'ee \`a des rappels sur les groupes de
r\'eflexions complexes imprimitifs, leurs repr\'esentations, et leurs
alg\`ebres de Hecke cyclotomiques. La construction alg\'ebrique
mentionn\'ee ci-dessus repose sur deux concepts: les repr\'esentations
sp\'eciales et l'induction tronqu\'ee. Nous les traitons aux
Sections~\ref{sect:special} et~\ref{sect:tronquee} respectivement. Nous
\'etablissons une compatibilit\'e entre l'induction tronqu\'ee et les
symboles \`a la Section~\ref{sect:Ge1n}, et une autre compatibilit\'e
entre les repr\'esentations sp\'eciales et les sous-groupes paraboliques
\`a la Section~\ref{sect:paraboliques}.

Enfin, \`a la Section~\ref{sect:springer}, nous d\'efinissons, de mani\`ere
alg\'ebrique, une classe de sous-groupes dits \emph{pseudoparaboliques} et
puis une classe de repr\'esentations dites \emph{de Springer}. (Pour les
groupes de Weyl, les re\-pr\'e\-sen\-ta\-tions de Springer sont celles
as\-so\-ci\'ees aux sys\-t\`e\-mes locaux triviaux par la correspondance de
Springer; elles sont donc en bijection avec les classes unipotentes).
Ensuite, nous calculons toutes les repr\'esentations de Springer de tous
les groupes imprimitifs sp\'etsiaux. Le Th\'e\-o\-r\`e\-me~\ref{thm:princ}
en donne un pa\-ra\-m\'e\-tra\-ge en termes des symboles dans le cas des
groupes non di\'edraux, et le Th\'eor\`eme~\ref{thm:princ-diedr} traite les
groupes di\'edraux.

%**************************************************************************
\section{Symboles et multipartitions}
\label{sect:symboles}
%**************************************************************************

Soient $d$ et $e$ deux entiers strictement positifs. Dans cette section, nous 
introduisons certains ensembles d'objets combinatoires (dont les \emph{symboles} 
et les \emph{multipartitions}) qui d\'ependent de $d$ et $e$. Dans la section 
suivante, nous rappellerons les liens entre ces objets et la th\'eorie des 
repr\'esentations du groupe de r\'eflexions complexes imprimitif $G(de,e,n)$ 
et de ses alg\`ebres de Hecke cyclotomiques. Pour cette raison, 
on dira toujours que 
nos objets combinatoires sont associ\'es au groupe $W = G(de,e,n)$, plut\^ot 
qu'aux entiers $d$ et $e$.

Un \emph{poids} pour $W$ est un \'el\'ement du quotient $
\Z^{de}/(1,\ldots,1)$, dont tout repr\'esent\-ant $(m_0, \ldots, m_{de-1})$
a la 
propri\'et\'e que $m_i = m_j$ si $i \equiv j \pmod d$. Par abus de langage, 
nous parlerons d'un \'el\'ement de $\Z^{de}$ comme s'il f\^ut un poids, au lieu 
de parler du poids dont cet \'el\'ement-l\`a est un repr\'esentant.

En particulier, le poids
\[
\bn(de,e) = (\underbrace{1,0,\ldots,0}_{\text{$d$ coordonn\'ees}},
\underbrace{1,0,\ldots,0}_{\text{$d$ coordonn\'ees}},\ldots
\underbrace{1,0,\ldots,0}_{\text{$d$ coordonn\'ees}})
\]
est appel\'e le \emph{poids sp\'etsial} pour $G(de,e,n)$. Les poids
\[
\bb = (1,0,\ldots, 0)
\qquad\text{et}\qquad
\bd = (0,\ldots,0)
\]
seront particuli\`erement utiles.

Soit $\Psi = (\Psi_0, \ldots, \Psi_{de-1})$ un $de$-uplet de suites finies 
croissantes d'entiers positifs:
\[
\Psi_i = (\Psi_i^{(0)} \le \cdots \le \Psi_i^{(m_i-1)}).
\]
Soient $r$ et $s$ deux entiers positifs, et supposons que $\Psi_i^{(0)} \ge s$ 
pour tout $i \ge 1$. On pose $\Psi' =
(\Psi'_0, \ldots \Psi'_{de-1})$, o\`u
\[
\Psi'_i =
\begin{cases} (0 \le \Psi_0^{(0)}+r \le \cdots \le \Psi_0^{(m_0-1)}+r) &
\text{si $i = 0$,} \\
(s \le \Psi_i^{(0)}+r \le \cdots \le \Psi_i^{(m_i-1)}+r) &
\text{si $i > 0$.}
\end{cases}
\]
On appelle $\Psi'$ le $(r,s)$-\emph{d\'ecal\'e} de $\Psi$. L'op\'eration de $
(r,s)$-d\'ecalage engendre une relation d'\'equivalence sur l'ensemble des $de$-uplets de suites finies croissantes d'entiers positifs. Une classe 
d'\'equivalence sous cette relation est appel\'ee un $(r,s)$-\emph{pr\'esymbole}. De plus, si $\bm$ est le poids $(m_0, \ldots, m_{de-1})$, 
on dit que $\Psi$ est \emph{de poids $\bm$}. (\'Evidemment, ce poids reste 
invariant sous d\'ecalage).

En particulier, le \emph{protosymbole} de type $(r,s)$ et de poids $\bm$ est le 
$(r,s)$-pr\'esymbole
\begin{align*}
\Phi &= \Phi^{r,s}(\bm) = (\Phi_0, \ldots, \Phi_{de-1}),\\
\intertext{o\`u}
\Phi_i &=
\begin{cases} (0, r, \ldots, (m_i-1)r) & \text{si $i = 0$,} \\
(s, r+s, \ldots, (m_i-1)r+s) & \text{si $i > 0$.}
\end{cases}
\end{align*}

Soit $\tilde\cP_n$ l'ensemble de $de$-uplets de partitions dont la somme totale \'egale $n$:
\[
\tilde\cP_n = \Bigg\{
( \underbrace{(0 \le \alpha_0^0 \le \cdots \le \alpha_0^{k_0})}_{\textstyle
\balpha_0},
\ldots, \underbrace{(0 \le \alpha_{de-1}^0 \le \cdots \le \alpha_{de-1}^{k_
{de-1}})}_{\textstyle\balpha_{de-1}}) \,\Bigg|\,
\sum_{i,j} \alpha_i^j = n \Bigg\}
\]
On appelle \emph{rotation} (\`a l'\'egard de $W$) l'application $\rot: \tilde\cP_n \to \tilde\cP_n$ d\'efinie par
\begin{equation}\label{eqn: rotation}
\rot(\balpha_0, \ldots, \balpha_{de-1}) = (\balpha_d, \balpha_{d+1}, \ldots, \balpha_{de-1}, \balpha_0, \ldots, \balpha_{d-1}).
\end{equation}
L'ensemble des \emph{multipartitions pour $W$}, ou \emph{$W$-multipartitions},  est l'ensemble des $\rot$-orbites sur $\tilde\cP_n$.

\'Evidemment, l'op\'eration de rotation est triviale dans le cas du groupe
$W = G(d,1,n)$, pour lequel une multipartition n'est autre qu'un $d$-uplet
de partitions. Par contre, les $\rot$-orbites sont en g\'en\'eral non
triviales pour $G(de,e,n)$. Si $\balpha$ est une
$G(de,e,n)$-multipartition, et si $\tbalpha \in \tilde\cP_n$ en est un
repr\'esentant, on dit que $\tbalpha$ est une $G(de,1,n)$-multipartition
\emph{au-dessus de} $\balpha$. En g\'en\'eral, lorsqu'on a besoin
d'\'ecrire une $G(de,e,n)$-multipartition explicitement, on \'ecrira
plut\^ot, par abus de notation, une $G(de,1,n)$-multipartition au-dessus de
celle-l\`a. Nous noterons $\cP(de,e,n)$ l'ensemble des
$G(de,e,n)$-multipartitions.

Soit $\balpha\in\cP(de,e,n)$, et soit $\tbalpha$ une
$G(de,1,n)$-multipartition au-dessus de $\balpha$. On note $s_e(\balpha)$
le cardinal du centralisateur de $\tbalpha$ dans le groupe cyclique
engendr\'e par $\rot$. (Il est clair que $s_e(\balpha)$ est ind\'ependant
du choix de $\tbalpha$). 

Ensuite, soit $\bm$ un poids. En ajoutant des ``$0$'' suppl\'ementaires si
n\'ecessaire, on peut consid\'erer $\tbalpha$ comme un $(0,0)$-pr\'esymbole
de poids $\bm$. Posons 
\[
\Lambda = \Lambda^{r,s}_\bm(\tbalpha) = \tbalpha + \Phi^{r,s}(\bm).
\]
Bien s\^ur, les divers $\tbalpha$ donnent lieu \`a divers $\Lambda$, et
$\rot$ induit une application (toujours appel\'ee rotation) sur l'ensemble
de tous les pr\'esymboles qui s'obtiennent de cette fa\c con.  Le
\emph{symbole de type $(r,s)$ et poids $\bm$} associ\'e \`a $\balpha$ est
l'ensemble de tous les pr\'esymboles obtenue par cette construction, ce qui
est une seule orbite par rotation. On note $Z^{r,s}_\bm$ l'ensemble des
symboles de type $(r,s)$ et de poids $\bm$.

Soit $\Lambda$ un symbole. Choisissons un repr\'esentant (\`a l'\'egard de
la rotation et du d\'ecalage) $(\Lambda_0, \ldots, \Lambda_{de-1})$ pour
$\Lambda$, ainsi qu'un repr\'esentant $(m_0, \ldots, m_{de-1}) \in \Z^{de}$
de son poids tel que $m_i$ \'egale le nombre de coefficients de
$\Lambda_i$. L'ensemble des ``positions'' dans $\Lambda$ est 
\[
\cS(\Lambda) = \{(i,j) \mid \text{$0 \le i < de$ et $0 \le j < m_i$}\}.
\]
Nous munissons cet ensemble d'un ordre total comme suit:
\[
(i,j) \prec (k,l)
\qquad\text{si}
\begin{cases}
\text{$j < l$}, &\text{ou} \\
\text{$j = l$, $k > 0$ et $i = 0$}, &\text{ou} \\
\text{$j = l$, $k > 0$ et $i > k$.}
\end{cases}
\]

Un symbole $\Lambda$ est \emph{distingu\'e} s'il poss\`ede un
repr\'esentant $(\Lambda_0, \ldots, \Lambda_{de-1})$ tel que
\begin{equation}\label{eqn:dist-defn}
\Lambda_i^{(j)} \le \Lambda_k^{(l)}
\qquad\text{si $(i,j) \prec (k,l)$}.
\end{equation}
(Autrement dit, un peut consid\'erer un repr\'esentant d'un symbole comme
une application $\cS(\Lambda) \to \N$; le symbole est distingu\'e si cette
application est croissante). Il est \`a noter que cette propri\'et\'e du
repr\'esentant est stable sous d\'ecalage mais non sous rotation en
g\'en\'eral. 

Deux symboles du m\^eme type et du m\^eme poids sont dits \emph{similaires} si 
tous deux poss\`edent des repr\'esentants ayant les m\^emes
coordonn\'ees avec les m\^emes multiplicit\'es. (Puisque les symboles sont
d\'efinis ici comme provenant des multipartitions, il n'est pas \'evident
que chaque classe de similitude contienne un symbole distingu\'e).

\begin{exm}
\newcommand{\vn}{\varnothing}
\newcommand{\st}[1]{%
\left(\begin{array}{@{}c@{}} #1\\ -\\ - \end{array}\right)}
\newcommand{\su}[4]{%
\left(\begin{array}{@{}c@{}c@{}c@{}}
#1&&#2\\ &#3&\\ &#4& \end{array}\right)}
\newcommand{\sv}[7]{%
\left(\begin{array}{@{}c@{}c@{}c@{}c@{}c@{}}
#1&&#2&&#3\\ &#4&&#5&\\ &#6&&#7& \end{array}\right)}
\newcommand{\sw}[6]{%
\left(\begin{array}{@{}c@{}c@{}c@{}c@{}c@{}c@{}c@{}}
0&&3&&6&&9\\ &#1&&#2&&#3& \\ &#4&&#5&&#6& \end{array}\right)}
Il y a 22 multipartitions pour $G(3,1,3)$:
{\tiny
\[
\begin{array}{llllll}
([3],\vn,\vn)   & ([2],\vn,[1])   & ([1],[1^2],\vn)
  & (\vn,[3],\vn)   & (\vn,[1^2],[1]) & (\vn,\vn,[1,2]) \\
([1,2],\vn,\vn) & ([1^2],[1],\vn) & ([1],[1],[1])
  & (\vn,[1,2],\vn) & (\vn,[1],[2])   & (\vn,\vn,[1^3]) \\
([1^3],\vn,\vn) & ([1^2],\vn,[1]) & ([1],\vn,[2])
  & (\vn,[1^3],\vn) & (\vn,[1],[1^2]) \\
([2],[1],\vn)   & ([1],[2],\vn)   & ([1],\vn,[1^2])
  & (\vn,[2],[1]) & (\vn,\vn,[3])
\end{array}
\]}
Les symboles correspondants de type $(3,1)$ et de poids $(1,0,0)$ sont:
{\tiny
\[
\begin{array}{cccccc}
\st3^*      & \su0512  & \sv0372514^*& \su0341^*   & \sv0362515^*& \sv0361426 \\
\su1511^*   & \su1421^*& \su0422^*   & \sv0362614^*& \su0323     & \sw147258 \\
\sv1471414^*& \su1412  & \su0413     & \sw258147^* & \sv0361525 \\
\su0521^*   & \su0431^*& \sv0371425  & \su0332^*   & \su0314
\end{array}
\]}
Les 13 symboles qui portent une \'etoile sont les symboles distingu\'es.

Un exemple d'une classe de similitude est:
{\tiny
\[
\left\{
\su0431, \su0341, \su0314\right\}
\]}
\end{exm}

\begin{exm}
\newcommand{\vn}{\varnothing}
\newcommand{\st}[3]{%
\left(\begin{array}{@{}c@{}} #1\\ #2\\ #3\end{array}\right)}
\newcommand{\su}[6]{%
\left(\begin{array}{@{}c@{\ }c@{}} #1&#2\\ #3&#4\\ #5&#6\end{array}\right)}
\newcommand{\sv}[9]{%
\left(\begin{array}{@{}c@{\ }c@{\ }c@{}}
#1&#2&#3\\ #4&#5&#6\\ #7&#8&#9\end{array}\right)}
Pour $G(3,3,3)$, l'op\'eration de rotation est non triviale. Voici un 
ensemble de repr\'esentants de ses 8 multipartitions:
{\tiny
\[
\begin{array}{llll}
([3],\vn,\vn)   & ([2],\vn,[1])   & ([1,2],\vn,\vn) & ([1^2],[1],\vn) \\
([1^3],\vn,\vn) & ([1^2],\vn,[1]) & ([2],[1],\vn)   & ([1],[1],[1])
\end{array}
\]}
Les symboles correspondants de type $(3,0)$ et de poids $(0,0,0)$ sont:
{\tiny
\[
\begin{array}{cccc}
\st300^*       & \st201    & \su150303^* & \su150403^* \\
\sv147036036^* & \su150304 & \st210^*    & \st111^*
\end{array}
\]}
Les \'etoiles d\'esignent toujours les symboles distingu\'es.
\end{exm}

%******************************************************************************
\section{Les groupes imprimitifs et leurs alg\`ebres de Hecke}
\label{sect:groupe}
%******************************************************************************

%------------------------------------------------------------------------------
\subsection{Les repr\'esentations irr\'eductibles du groupe $G(de,e,n)$}
\label{subsection: repr}
%------------------------------------------------------------------------------

Rappelons que le groupe complexe imprimitif $G(e,1,n)$ est le groupe 
lin\'eaire complexe sur $V = \bigoplus_{j=1}^n\C e_j$ form\'e des matrices 
monomiales dont les coefficients non nuls appartiennent \`a $\{\zeta_e^j:0\le j
\le e-1\}$, o\`u $\zeta_e$ est une racine primitive $e$-i\`eme de l'unit\'e. Le
groupe $G(e,1,n)$ est donc le produit semi-direct de son sous-groupe de matrices 
diagonales avec le sous-groupe de matrices de permutations, \ie $G(e,1,n)=(\Z/e
\Z)^n\rtimes\fS_n$. Dans cette repr\'esentation, le groupe $G(e,1,n)$ est 
engendr\'e par la r\'eflexion $t$ qui envoie $e_1$ sur $\zeta_e e_1$ et laisse 
fixes $e_2$, $\ldots$, $e_n$ et par les matrices de permutations $s_i$ ($1\le i
\le n-1$) correspondant aux transpositions $(i,i+1)$. 

Soit $\gamma_e\colon G(e,1,n)\to\C$ le caract\`ere lin\'eaire d\'efini par 
$\gamma_e(t):=\zeta_e$ et $\gamma_e(s_i):=1$ pour $1\le i\le n-1$.
 
Soit $\balpha=(\balpha_0, \balpha_1, \ldots,\balpha_{e-1})$ un $e$-uplet 
de partitions de $n$. Pour tout entier $i$ tel que $0\le i\le e-1$, nous notons 
$n_i$ la \emph{somme de la partition}
$\balpha_i$ (\ie $n_i:=\sum_{j=0}^{k_i}\balpha_i^j$). Les repr\'esentations du 
groupe sym\'etrique $\fS_{n_i}$ peuvent
\^etre consid\'er\'ees comme des repr\'esentations du groupe $G(e,1,n_i)$, via 
la projection naturelle de $G(e,1,n_i)$
sur $\fS_{n_i}$. Les classes d'isomorphie des repr\'esentations irr\'eductibles 
de $\fS_{n_i}$ sont param\'etr\'ees par 
les partitions de $n_i$ et nous noterons 
$E_{\balpha_i}$ une repr\'esentation ir\-r\'e\-duc\-ti\-ble de $\fS_{n_i}$ 
correspondant \`a la partition $\balpha_i$ de $n_i$. Nous posons $\bn_e:=
(n_0,n_1,\ldots,n_{e-1})$ et
\[
G(e,1,\bn_e):=G(e,1,n_0)\times\cdots\times G(e,1,n_{e-1}).
\] 
La formule
\[
E_\balpha:=
\Ind^{G(e,1,n)}_{G(e,1,\bn_e)}
\left(E_{\balpha_0}\otimes(E_{\balpha_1}\otimes\gamma_e)\otimes
\cdots\otimes(E_{\balpha_{e-1}}\otimes\gamma_e^{e-1})\right)
\]
d\'efinit donc une repr\'esentation du groupe $G(e,1,n)$. La repr\'esentation 
$E_\balpha$ est irr\'eductible,
$E_\balpha\not\simeq E_\bbeta$ si 
$\balpha\ne\bbeta$, et les (classes d'isomorphie) des $E_\balpha$ d\'ecrivent 
toutes les (classes d'isomorphie) de repr\'esentations 
irr\'eductibles de $G(e,1,n)$.

\smallskip Le groupe $G(de,e,n)$ est un sous-groupe d'indice $e$ de $G(de,1,n)$, 
noyau du caract\`ere lin\'eaire
$\gamma_{de}^d$. Nous allons rappeler la description de ses caract\`eres 
ir\-r\'eductibles en fonction de ceux du groupes
$G(de,1,n)$.

Remarquons que $s_e(\balpha)$ (voir la Section~\ref{sect:symboles}) divise $n$. La restriction de $E_\balpha$ \`a $G(de,e,n)$ est 
somme de $s_e(\balpha)$ repr\'esentations
irr\'eductibles distinctes, nous les notons $E_{\balpha,1}$, $\ldots$, $E_
{\balpha,s_e(\balpha)}$ et toute
repr\'esentation irr\'eductible de $G(de,e,n)$ intervient dans la restriction 
d'une repr\'esentation $E_\balpha$ pour
$\balpha$ un $de$-uplet de partitions de $n$. 
Plus pr\'ecis\'ement, nous posons 
\[
\sigma:=(s_1s_2\cdots s_{n-1})^{n/s_e(\balpha)}\quad \text{et}\quad G(de,e,\bn_
{de}):=G(de,1,\bn_{de})\cap G(de,e,n).
\]
La restriction
\`a $G(de,e,\bn_{de})$ de la repr\'esentation 
\[
E_{\balpha_0}\otimes(E_{\balpha_1}\otimes\gamma_{de})\otimes
\cdots\otimes(E_{\balpha_{de-1}}\otimes\gamma_{de}^{de-1})
\]
de $G(de,1,\bn_{de})$ est invariante par $\sigma$ et s'\'etend au produit semi-
direct
$G(de,e,\bn_{de})\rtimes\langle\sigma\rangle$. Les induites \`a $G(de,e,n)$ des 
diverses extensions d\'ecrivent
l'ensemble des composantes irr\'eductibles de la restriction de $E_\balpha$ \`a 
$G(de,e,n)$.

%------------------------------------------------------------------------------
\subsection{Polyn\^ome de Poincar\'e et degr\'es fant\^omes} 
\label{subsection: fantomes}
%------------------------------------------------------------------------------

Soit $W\subset\GL(V)$ un groupe de r\'eflexions complexes et soit $S(V)$ 
l'al\-g\`e\-bre sy\-m\'e\-tri\-que 
de $V$. Nous notons $n$ la dimension de $V$.
L'al\-g\`e\-bre des in\-va\-riants $S(V)^W$ de $W$ dans $S(V)$ 
est une al\-g\`e\-bre de po\-ly\-n\^o\-mes sur $n$ \'e\-l\'e\-ments 
ho\-mo\-g\`e\-nes al\-g\'e\-bri\-que\-ment in\-d\'e\-pen\-dants de de\-gr\'es 
respectifs no\-t\'es $d_1$, $\ldots$, $d_n$ (\cite{Ch}). 
Le nombre $N^*$ de r\'eflexions de $W$ est \'egal \`a $\sum_{i=1}^n(d_i-1)$. 

Le \emph{polyn\^ome de Poincar\'e} $P_W$ de $W$ est donn\'e par la formule
\[
(X-1)^n\cdot P_W(X)=\left(\frac{1}{|W|}
\sum_{w\in W}\frac{\det_V(w)}{\det_V(X-w)}\right)^{-1} =\prod_{i=1}^n(X^
{d_j}-1),
\] 
o\`u $\det_V$ d\'esigne le d\'eterminant sur $V$. 
Pour $W=G(de,e,n)$, on obtient
\[
P_{G(de,e,n)}(X)=\frac{X^{dn}-1}{X-1}\cdot
\prod_{i=1}^{n-1}\frac{X^{dei}-1}{X^i-1}.
\]
En particulier:
\[
P_{G(e,1,n)}(X)=\prod_{i=1}^n\frac{X^{ei}-1}{X-1},\quad P_{G(e,e,n)}(X)=
\frac{X^{n}-1}{X-1}\cdot\prod_{i=1}^{n-1}\frac{X^{ei}-1}{X-1}.
\]

\smallskip

Soit $S(V)^W_+$ l'id\'eal de $S(V)^W$ form\'e des \'el\'ements de degr\'es
strictement positifs. Nous notons $S(V)_W:=S(V)/(S(V)^W_+\cdot S(V))$
l'\emph{alg\`ebre coinvariante} de $(W,V)$. En tant que $W$-module,
$S(V)_W$ est isomorphe \`a la repr\'esentation r\'eguli\`ere de $W$.
Soit $S(V)_W=\bigoplus_{j=0}^{N^*}S(V)_W^j$ la d\'ecomposition
de $S(V)_W$ en ses composantes gradu\'ees.
 
Le \emph{degr\'e fant\^ome} d'une repr\'esentation irr\'eductible $E$ de
$W$, not\'e $R_E(X)$, est le polyn\^ome dans $\Z[X]$ d\'efini par 
\begin{equation} \label{eqdef: RE}
R_E(X):=\sum_j m_j(E)X^j,
\end{equation}
o\`u $m_j(E)$ d\'esigne la multiplicit\'e avec laquelle $E$ appara\^\i t
dans le $W$-module $S(V)_W^j$.
On a (voir par exemple \cite{Sp}):
\begin{equation} \label{eq: RE}
R_E(X)=(X-1)^n\cdot P_W(X)\cdot 
\sum_{w\in W}\frac{\det_V(w)\Tr(w,E)}{\det_V(X-w)}.
\end{equation}

\smallskip

Pour toute partie finie $A$ de $\N$, nous d\'efinissons les polyn\^omes 
suivants:
\[
\Delta(A,X):=\prod_{\substack{a,b \in A\\ b<a}}(X^a-X^b),
\]
\[
\Theta(A,X):=\prod_{a\in A}\,\prod_{l=1}^a(X^l-1).
\]
Soit $\balpha=(\balpha_0, \balpha_1, \ldots, \balpha_{de-1})$ un $de$-uplet de 
partitions de $n$. Pour $0\le i\le de-1$, nous \'ecrivons 
$\balpha_i=(0\le\alpha_i^0\le\cdots\le\alpha_i^{k_i})$, nous notons
$n_i(\balpha)$ la somme de la partition $\balpha_i$, posons
\[c_i(\balpha):=\sum_{l=0}^{k_i}\binom{l}{2},\]
et d\'efinissons des parties finies 
\[
A_i(\balpha):=\left\{\bar\balpha_i^0,\bar\balpha_i^1,\ldots,\bar\balpha_i^{k_i}
\right\}
\]
de $\N$, o\`u 
\[
\bar\balpha_i^j :=\balpha_i^j+j, \quad\text{pour $0\le j\le k_i$.}
\]
Pour tout multiple $m$ de $de$, nous poserons $\balpha_{i+m}:=\alpha_i$, $c_{i
+m}(\balpha)=c_i(\balpha)$ et
$A_{i+m}(\balpha)=A_i(\balpha)$.

\smallskip

D'apr\`es \cite[Remarque~2.10]{M1}, le degr\'e fant\^ome $R_{E_\balpha}(X)$ de 
la 
repr\'esentation irr\'e\-ductible $E_\balpha$ du groupe $G(de,1,n)$
s'\'ecrit:
\begin{equation}
\label{eqn: fantome} R_{E_\balpha}(X)=\prod_{h=1}^n(X^{deh}-1)\cdot
\prod_{i=0}^{de-1}\frac{\Delta(A_i(\balpha),X^{de})\cdot X^{in_i}}{\Theta(A_i
(\balpha),X^{de})\cdot X^{de\cdot
c(\balpha_i)}}.
\end{equation}
D'apr\`es \cite[p.~806]{M1}, les repr\'esentations irr\'eductibles 
$E_{\balpha,l}$, pour $l\in\{1,\ldots, s_e(\balpha)\}$, ont toutes le m\^eme 
degr\'e fant\^ome $R_{E_{\balpha,l}}(X)=:R_{E_{\balpha},d}$, lequel s'\'ecrit
\[
R_{E_{\balpha},d}(X)=\frac{X^{nd}-1}{X^{nde}-1}\cdot\frac{1}{s_e(\balpha)}\cdot
\sum_{j=0}^{e-1}R_{\rot^{j}(\balpha)}(X), 
\]
o\`u $\rot$ est d\'efinie par~(\eqref{eqn: rotation}). Puisque 
\[
\rot^j(\balpha)=(\balpha_{jd},\balpha_{jd+1},\ldots,\balpha_{de-1},\balpha_0,
\balpha_1,\ldots,\balpha_{jd-1}),
\]
nous obtenons 
\[
(\rot^j\balpha)_i=\balpha_{i+jd},\;\;\text{pour $0\le i\le de-1$ et $0\le j\le 
e-1$.}
\]
Le polyn\^ome $R_{E_{\balpha},d}(X)$ admet donc l'expression 
\begin{multline*} \frac{X^{nd}-1}{X^{nde}-1}\cdot\frac{1}{s_e(\balpha)}\cdot
\sum_{j=0}^{e-1}\left(
\prod_{h=1}^n(X^{deh}-1)\cdot
\prod_{i=0}^{de-1}\frac{\Delta(A_{i+jd}(\balpha),X^{de})\cdot
X^{in_{i+jd}(\balpha)}}{\Theta(A_{i+jd}(\balpha),X^{de})\cdot 
X^{de\cdot c_{i+jd}(\balpha)}}\right)\\
=\frac{(X^{nd}-1)}{s_e(\balpha)}\cdot\prod_{h=1}^{n-1}(X^{deh}-1)\cdot
\sum_{j=0}^{e-1}\prod_{i=0}^{de-1}\frac{\Delta(A_{i+jd}(\balpha),X^{de})\cdot
X^{in_{i+jd}(\balpha)}}{\Theta(A_{i+jd}(\balpha),X^{de})\cdot 
X^{de\cdot c_{i+jd}(\balpha)}}.
\end{multline*}
L'expression
\[
\prod_{i=0}^{de-1}\frac{\Delta(A_{i+jd}(\balpha),X^{de})} {\Theta(A_{i+jd}
(\balpha),X^{de})X^{de\cdot
c_{i+jd}(\balpha)}}
\]
\'etant ind\'ependante du choix de $j\in\{0,\ldots,e-1\}$, le polyn\^ome $R_{E_
{\balpha},d}(X)$ s'\'ecrit encore
\[
\frac{(X^{nd}-1)}{s_e(\balpha)}\cdot\prod_{h=1}^{n-1}(X^{deh}-1)\cdot
\prod_{i=0}^{de-1}\frac{\Delta(A_{i}(\balpha),X^{de})} {\Theta(A_{i}(\balpha),X^
{de})X^{de\cdot c_{i}(\balpha)}}\cdot
\sum_{j=0}^{e-1}\prod_{i=0}^{de-1}X^{in_{i+jd}(\balpha)}.
\]
En utilisant les \'egalit\'es
\[
n_{i+jd}(\balpha)=\sum_{l=0}^{k_{i+jd}}(\bar\balpha_{i+jd}^l-l)= -\frac{k_{i+jd}
(k_{i+jd}+1)}{2}+\sum_{a\in
A_{i+jd}(\balpha)}a,
\]
nous voyons que $R_{E_{\balpha},d}(X)$ est \'egal \`a
\begin{multline*}
\frac{(X^{nd}-1)}{s_e(\balpha)}\cdot\prod_{h=1}^{n-1}(X^{deh}-1)\cdot
\prod_{i=0}^{de-1}\frac{\Delta(A_{i}(\balpha),X^{de})} {\Theta(A_{i}(\balpha),X^
{de}) X^{de\cdot c_i(\balpha)}} \\
\cdot
\sum_{j=0}^{e-1}\prod_{a\in A_{i+jd}(\balpha)}X^{i(a-\frac{k_{i+jd}(k_{i+jd}+1)}
{2})}.
\end{multline*}

%------------------------------------------------------------------------------
\subsection{Alg\`ebres de Hecke cyclotomiques et familles de caract\`eres}
%------------------------------------------------------------------------------
Soit $W$ un groupe de r\'eflexions complexes irr\'eductible fini et soit $\cD$ 
le diagramme qui lui est associ\'e dans \cite{BMR}. Ceci d\'efinit une
pr\'esentation de $W$ sur un ensemble de g\'en\'erateurs $S$, avec des
``relations d'ordre'' $s^{d_s}=1$ pour $s\in S$, ainsi que des relations
homog\`enes, appel\'ees ``relations de tresses''. Le ``groupe de tresses''
$B=B(W)$ associ\'e \`a $W$ est par d\'efinition le groupe engendr\'e par
un ensemble $\{\bs\,:\,s\in S\}$ en bijection $\bs\leftrightarrow s$ avec
$S$, satisfaisant aux relations de tresses de $\cD$. Soit
$\bu=\{u_{s,i}\,:\,s\in S, 0\le i\le d_s-1\}$ un ensemble
de nombres transcendants sur $\Z$ tels que $u_{s,i}=u_{t,i}$ si $s$ et $t$
sont conjugu\'es dans $W$. L'\emph{alg\`ebre de Hecke g\'en\'erique}
$\cH(W,\bu)$ de $W$ de param\`etre $\bu$ est d\'efinie comme le quotient 
$$\cH(W,\bu):=\Z[\bu,\bu^{-1}]B/I, \quad\text{ avec
$I=\left(\prod_{i=0}^{d_s-1}(\bs-u_{s,i})\,:\,s\in S\right)$}$$
de l'alg\`ebre de groupe de $B$ sur $\Z[\bu,\bu^{-1}]$ par l'id\'eal
$I$ engendr\'e par certaines ``relations d'ordre d\'eform\'ees''.

Soient $\bmu_\infty$ le sous-groupe des racines de l'unit\'e de $\C$ et 
$K$ un sous-corps du corps $\Q(\bmu_\infty)$ de degr\'e fini sur $\Q$. On
note $\Z_K$ l'anneau des entiers de $K$ (c'est un anneau de Dedekind) et
$\bmu(K)$ le groupe des racines de l'unit\'e de $K$. Soit $\zeta$
un \'el\'ement de $\bmu(K)$.
Pour $s\in S$ et $0\le i\le d_s-1$, nous supposons donn\'es des entiers
relatifs $n_{s,i}\in\Z$. Nous posons 
$m_{s,i}:=n_{s,i}/|\bmu(K)|$ et $\bm_s:=(m_{s,0},m_{s,1},\ldots,m_{s,d_s-1})$, et
nous notons $\bm$ l'ensemble $\{\bm_s\,:\,s\in S\}$.
L'\emph{alg\`ebre de Hecke $\zeta$-cyclotomique} $\cH_\zeta^\bm(W)$ de $W$ est la 
$\Z_K[q,q^{-1}]$-alg\`ebre
obtenue \`a partir de $\cH(W,\bu)$ au moyen de la sp\'ecialisation
$\phi\colon\Z_K[\bu,\bu^{-1}]\to\Z_K[q,q^{-1}]$ d\'efinie par
$$\phi\colon u_{s,i}\mapsto \zeta_{d_s}^i(\zeta^{-1}q)^{m_{s,i}}.$$
  
Le \emph{degr\'e g\'en\'erique} d'une repr\'esentation irr\'eductible $E$
de $W$ est le polyn\^ome
$$D_E=P_W/c_E$$
quotient du polyn\^ome de Poincar\'e de $W$ par l'el\'ement de Schur $c_E$ de
$E$.

L'alg\`ebre de Hecke g\'en\'erique $\cH(e,1,n)$ associ\'ee au groupe $G(e,1,n)$ 
est l'alg\`ebre engendr\'ee sur l'anneau des polyn\^omes
de Laurent en $e+2$ ind\'etermin\'ees:
$$\Z[u_0,u_1,u_0^{-1},u_1^{-1},v_0,v_1,\ldots,v_{e-1},v_0^{-1},v_1^{-1},\ldots,v_{e-1}^{-1}]$$
par des \'el\'ements $s_1$, $s_2$, $\ldots$, $s_{n-1}$, $t$ satisfaisant
les relations de tresses
\[s_js_{j+1}s_j=s_{j+1}s_js_{j+1} \quad\text{et}\quad
s_{n-1}ts_{n-1}t=ts_{n-1}ts_{n-1}\] et les relations relations d'ordre
d\'eform\'ees
$$(s_j-u_0)(s_j-u_1)=(t-v_0)(t-v_1)\cdots(t-v_{e-1})=0.$$ 
L'\emph{alg\`ebre sp\'etsiale} de $G(e,1,n)$ est l'alg\`ebre de Hecke $1$-cyclotomique,
obtenue par la sp\'e\-cia\-li\-sa\-tion $1$-cyclotomique:
$$u_0\mapsto q,\quad u_1\mapsto -1,\quad
v_0\mapsto q \quad\text{et}\quad
v_i\mapsto \zeta_e^i\quad\text{pour $1\le i\le e-1$,}$$
\ie l'alg\`ebre $\cH^\bm_1(e,1,n)=\cH^\bm_1(G(e,1,n))$, avec 
$\bm=\{\bm_{s_1},\ldots,\bm_{s_{n-1}},\bm_t\}$, o\`u $\bm_{s_j}=(1,0)$ pour 
$1\le j\le n-1$ et $\bm_t=(1,0,\ldots,0)$.

\smallskip

L'alg\`ebre de Hecke g\'en\'erique $\cH(e,e,n)$ associ\'ee au groupe
$G(e,e,n)$ est:
\begin{itemize}
\item
si $n>2$ ou $n=2$ et $e$ impair, l'al\-g\`e\-bre en\-gen\-dr\'ee sur l'an\-neau 
$$\Z[u_0,u_1,u_0^{-1},u_1^{-1}]$$ par des \'e\-l\'e\-ments $s_1$, 
$s_2$, $\ldots$, $s_{n-1}$, $s_{n-1}'$ satisfaisant les relations 
$$s_js_{j+1}s_j=s_{j+1}s_js_{j+1}\;\;\text{($1\le j\le n-2$)}, 
\quad s'_{n-1}s_{n-2}s_{n-1}'=s_{n-2}s'_{n-1}s_{n-2},$$ 
$$s_{n-2}s_{n-1}'s_{n-1}s_{n-2}s'_{n-1}s_{n-1}=s_{n-1}'s_{n-1}s_{n-2}s'_{n-1}
s_{n-1}s_{n-2},$$
$$
\underbrace{s_{n-1}s_{n-1}'s_{n-1}s_{n-1}'s_{n-1}s_{n-1}'\cdots}_{\text{$e$ facteurs}}
=\underbrace{s_{n-1}'s_{n-1}s_{n-1}'s_{n-1}s'_{n-1}s_{n-1}\cdots}_{\text{$e$ facteurs}}$$
$$\text{et}\quad (s_{n-1}'-u_0)(s'_{n-1}-u_1)=(s_j-u_0)(s_j-u_1)=0,\quad
\text{ pour $1\le i\le n-1$};$$
\item
si $n=2$ et $e$ pair, l'al\-g\`e\-bre en\-gen\-dr\'ee sur l'an\-neau
$$\Z[u_0,u_1,v_0,v_1,u_0^{-1},u_1^{-1},v_0^{-1},v_1^{-1}]$$ 
par des \'e\-l\'e\-ments $s_1$, $s'_1$ satisfaisant les relations 
$$\underbrace{s_1s_1's_1s_1's_1s_1'\cdots}_{\text{$e$ facteurs}}
=\underbrace{s_1's_1s_1's_1s'_1s_1\cdots}_{\text{$e$ facteurs}}$$
$$\text{et}\quad (s_1'-u_0)(s'_1-u_1)=(s_1-v_0)(s_1-v_1)=0.$$
\end{itemize}
L'\emph{alg\`ebre sp\'etsiale} de $G(e,e,n)$ est l'alg\`ebre de Hecke $1$-cyclotomique,
obtenue par la sp\'e\-cia\-li\-sa\-tion $1$-cyclotomique:
$$u_0\mapsto q,\quad u_1\mapsto -1,\quad v_0\mapsto q\quad\text{et}\quad
v_1\mapsto -1,$$
\ie l'alg\`ebre $\cH^\bm_1(e,e,n)=\cH^\bm_1(G(e,e,n))$, avec 
$\bm=\{\bm_{s'_1},\bm_{s_1},\ldots,
\bm_{s_{n-1}}\}$, o\`u $\bm_{s'_1}=\bm_{s_1}=\cdots=\bm_{s_{n-1}}=(1,0)$.

%Expliquer le lien entre les poids et les sp\'ecialisations $\zeta$-cyclotomiques
%de l'alg\`ebre de Hecke g\'en\'erique. 

\smallskip

Lusztig a construit une partition des caract\`eres irr\'eductibles d'un
groupe de Coxeter fini $W$ en \emph{familles} \`a l'aide de la th\'eorie des
\emph{cellules}. Cette partition appara\^\i t naturellement dans le param\'etrage
de Lusztig des caract\`eres unipotents d'un groupe r\'eductif sur un corps
fini. Pour le moment il n'existe pas de d\'efinition de cellules pour les
groupes de r\'eflexions complexes et l'on ne peut donc pas utiliser
l'approche de Lusztig pour d\'efinir les familles de caract\`eres. Dans
\cite{R}, Rouquier a d\'ecrit une approche diff\'erente dans laquelle les
familles sont d\'efinie comme les blocs de d'alg\`ebre de Hecke-Iwahori de
$W$ sur un certain anneau $\cO(x)$, d\'efini comme suit:
$$\cO(x):=\Z_K[x,x^{-1},(x^m-1)^{-1}_{m\ge 1}].$$

Nous supposons dor\'enavant donn\'e un groupe de r\'eflexions complexes $W$
imprimitif et nous fixons une bijection de l'ensemble des
caract\`eres irr\'eductibles de $W$ sur celui des
caract\`eres irr\'eductibles de $\cH^\bm_\zeta(W)$ et identifions
ces deux ensembles via la bijection. Les \emph{familles} 
des caract\`eres irr\'eductibles sont alors d\'efinies (voir
\cite[Definition~2.4]{MR}) comme les blocs de $\cO(x)\cH^\bm_\zeta(W)$. 

Brou\'e et Kim ont d\'emontr\'e (voir \cite[Th\'eor\`eme~3.17]{BK}) que les 
familles de caract\`eres de $G(e,1,n)$ \`a l'\'egard de l'alg\`ebre
sp\'etsiale sont donn\'ees par les classes de 
similitude des symboles de type $(1,0)$ et de poids $\bb$.

Afin de d\'ecrire les familles de caract\`eres de $G(e,e,n)$ \`a l'\'egard
de l'alg\`ebre sp\'etsiale, introduisons
la notion de $e$-partition \emph{b\'egayante}: une $e$-partition $\balpha=
(\balpha_0,\ldots,\balpha_{e-1})$ est dite b\'egayante si
$\balpha_0=\cdots=\balpha_{e-1}$.
\`A toute $e$-partition b\'egayante de somme $n$ correspond $e$ familles
de caract\`eres de $G(e,e,n)$, chacune r\'eduite \`a un singleton. Les
autres familles de caract\`eres de $G(e,e,n)$ sont donn\'ees par les
classes de similitude des symboles associ\'es \`a des $e$-partition
non  b\'egayantes, de type $(1,0)$ et de poids $\bd$ (voir
\cite[Th\'eor\`eme 4.3]{BK}).  
%------------------------------------------------------------------------------
\subsection{Exemples}
%------------------------------------------------------------------------------

Les familles de caract\`eres pour les groupes de Weyl classiques $G(2,e,n)$ ($e = 
1,2$) correspondent aux classes de
similitude dans $Z^{1,0}_{\bn(2,e)}$, et en particulier, les caract\`eres 
sp\'eciaux correspondent aux symboles distingu\'es.

Les classes unipotentes pour $B_n$ (resp. $C_n$, $D_n$) correspondent aux 
classes de similitude dans
$Z^{2,0}_{\bn(2,1)}$ (resp.~$Z^{2,1}_{\bn(2,1)}$, $Z^{2,0}_{\bn(2,2)}$). En 
particulier, les caract\`eres associ\'es aux
syst\`emes locaux triviaux par la correspondence de Springer correspondent aux 
symboles distingu\'es.

Les calculs de Brou\'e--Kim \cite{BK} et de Kim \cite{K} montrent que les 
familles de caract\`eres (\`a 
l'\'egard d'une alg\`ebre cyclotomique) de
$G(de,e,n)$ sont en bijection avec les classes de similitude de symboles de 
type~$(1,0)$ et de poids convenable. Dans le
cas de l'alg\`ebre sp\'etsiale d'un groupe sp\'etsial, on sait aussi, d'apr\`es 
Malle, que les caract\`eres sp\'eciaux
correspondent aux symboles distingu\'es de type~$(1,0)$.

\begin{exm} \label{exm: cyclique}
Consid\'erons le cas tr\`es simple du groupe $G(e,1,1)$ (lequel est un
groupe cyclique d'ordre $e$). Nous savons d\'ej\`a par
\cite[Proposition~2.10(2)]{BK} qu'il y a deux familles de
repr\'esentations irr\'eductibles de $G(e,1,1)$, celle r\'eduite \`a la
repr\'esentation triviale et celle form\'ee des autres repr\'esentations
irr\'eductibles. Ce r\'esultat se r\'einterpr\`ete en termes de symboles
de la mani\`ere suivante.
On a $\Phi^{1,0}(\bb)=\hbox{\tiny$\left(\begin{array}{ll}0&1\cr 0&\cr\vdots\cr
0\end{array}\right)$}$.
Il y a $e$ multipartitions pour $G(e,1,1)$. Les symboles correspondants
de type $(1,0)$ et de poids $\bb$ sont
$$\Lambda_0=
\hbox{\tiny$\left(\begin{array}{l}1\cr-\cr\vdots\cr-\end{array}\right)$}
\quad\text{et}\quad
\Lambda_i=\hbox{\tiny$\left(\begin{array}{ll}0&1\cr 0&\cr\vdots&\cr 0\cr 1\cr 0\cr\vdots\cr
0\end{array}\right)
\begin{array}{l}
\leftarrow \text{$0$-\`eme ligne} \\
\text{ }\\ \qquad\vdots \\ \text{ }\\
\leftarrow \text{$i$-\`eme ligne ($1 \le i \le e-1$).} \\
\text{ }\\ \qquad\vdots \\ \text{ }
\end{array}$}$$
Les symboles sp\'eciaux sont $\Lambda_0$ et $\Lambda_1$. Les classes de
similitude sont au nombre de deux: $$\{\Lambda_0\}\quad\text{et}\quad
\{\Lambda_i\,:\,1\le i\le e-1\}.$$ 
%On a $a([1,-,\ldots,-])=0$ et $a([-,\ldots,-,i,-,\dots,-])=1$ pour $i\ne 0$.
D'autre part, la repr\'esentation induite 
$\Ind_{\{1\}}^{G(e,1,1)}(1)$, \'etant \'egale \`a la repr\'esentation
r\'eguli\`ere du groupe $G(e,1,1)$,  est somme de la repr\'esentation
triviale et de la repr\'esentation
$\gamma_e\oplus\gamma_e^2\oplus\cdots\oplus\gamma_e^{e-1}$. Ces deux
repr\'esentations constituent donc les repr\'esentations ``constructibles'' du
groupe $G(e,1,1)$.
\end{exm}

%**************************************************************************
\section{Repr\'esentations sp\'eciales}
\label{sect:special}
%**************************************************************************

%--------------------------------------------------------------------------
\subsection{Les fonctions $a$ et $b$}
\label{sect:ab}
%--------------------------------------------------------------------------

Soit $W$ un groupe de r\'eflexions complexes fini, et choisissons une alg\`ebre
de Hecke cyclotomique pour $W$. Si $E$ est une 
repr\'esentation irr\'eductible de $W$, on peut consid\'erer la multiplicit\'e 
de la racine en $q = 0$ de son degr\'e g\'en\'erique et de son degr\'e
fant\^ome. Ces deux entiers, qui sont appel\'es ``$a(E)$'' et
``$b(E)$'', respectivement, dans la litt\'erature, jouent un r\^ole tr\`es
important dans la suite. La repr\'esentation $E$ est dite \emph{sp\'eciale} si
$a(E)=b(E)$.

Dans le cas o\`u $W$ est imprimitif, ces fonctions ne d\'ependent que de la
multipartition associ\'ee \`a $E$. Nous \'ecrirons donc ``$a(\balpha)$'' et
``$b(\balpha)$'' au lieu de ``$a(E_{\balpha,l})$'' et
``$b(E_{\balpha,l})$''. 

Les fonctions $a$ et $b$ permettent aussi de d\'efinir la notion, introduite
dans \cite{M2}, de groupe de r\'eflexions complexes fini \emph{sp\'etsial}:
un groupe de r\'eflexions complexes fini $W$ est dit sp\'etsial si l'on a
\begin{equation} \label{spetsialite}
a(E)\le b(E)\quad\text{pour toute repr\'esentation irr\'eductible $E$ de
$W$.}\end{equation}
D'apr\`es \cite[Proposition~8.1]{M3}, la propri\'et\'e~(\ref{spetsialite}) 
est \'equivalente au fait que
pour toute repr\'esentation irr\'eductible $E$ de
$W$, il existe une repr\'esentation irr\'eductible sp\'eciale $E_0$
telle que $a(E)=a(E_0)$.
La Proposition~8.1 de \cite{M3} fournit d'autres caract\'erisations des
groupes de r\'eflexions complexes finis sp\'etsiaux.

L'ensemble des groupes de r\'eflexions complexes finis sp\'etsiaux contient en
particulier tous les groupes de r\'eflexions complexes finis qui peuvent \^etre 
d\'efinis sur le corps des nombres r\'eels.

Les groupes de r\'eflexions complexes finis sp\'etsiaux imprimitifs 
irr\'eductibles sont les groupes $\fS_n$, $G(e,1,n)$ et $G(e,e,n)$ (voir 
\cite[preuve de la Proposition~8.1]{M3}).

\smallskip
 
Nous rappelons maintenant la formule obtenue par Malle~\cite{M1} pour
$b(\balpha)$. Il r\'esulte de l'expression obtenue pour $R_{E_
{\balpha},d}(X)$ que $b(\balpha)$ est
\'egal \`a
\begin{equation}
\label{eqn:b} ed\sum_{i=0}^{de-1}\left(\sum_{\substack{a,b\in
A_i(\balpha)\\ b<a}}b-c_i
\right) -\min_{0\le j\le e-1}
\sum_{i=0}^{de-1}i\left(\frac{k_{i+jd}(k_{i+jd}+1)}{2}-
\sum_{a\in A_{i+jd}(\balpha)}a\right).
\end{equation}
Lorsque $W=G(e,1,n)$ et que la suite $(k_0,k_1,\ldots,k_{e-1})$ est telle 
que $k_0=k+1$ et $k_i=k$ si $1\le i\le e-1$, avec $k\in\N$, la 
formule~(\ref{eqn:b}) donne
\begin{equation}
\label{eqn: bGe1n} b(\balpha)=e\sum_{i=0}^{e-1}\sum_{\substack{a,b\in
A_i(\balpha)\\ 
b<a}}b+\sum_{i=0}^{e-1}i\sum_{a\in
A_{i}(\balpha)}a -\sum_{l=0}^{k-1}\binom{el+1}{2}.
\end{equation}
Lorsque $W=G(e,e,n)$ et que la suite $(k_0,k_1,\ldots,k_{e-1})$ est telle 
que $k_i=k$ pour tout $i\in\{0,\ldots,e-1\}$, avec $k\in\N$, la formule~(\ref
{eqn:b}) donne
\begin{equation}
\label{eqn: bGeen} b(\balpha)=e\sum_{i=0}^{e-1}\sum_{\substack{a,b\in
A_i(\balpha)\\ 
b<a}}b +\min_{0\le j\le e-1}
\sum_{i=0}^{e-1}i
\sum_{a\in A_{i+jd}(\balpha)}a -\sum_{l=0}^{k-1}\binom{el}{2}.\end{equation}

%--------------------------------------------------------------------------
\subsection{Formules combinatoires pour $a$ et $b$}
%--------------------------------------------------------------------------

Dans cette section, nous introduisons et \'etudions certaines fonctions
combinatoires d\'efinies sur l'ensemble de symboles pour un groupe
imprimitif. Les quelques derniers \'enonc\'es traitent la relation entre
ces fonctions combinatoires et les fonctions $a$ et $b$ de la section
pr\'ec\'edente.

Soit $\balpha \in \cP(e,1,n)$, et soit $\Lambda$ son symbole de
type $(r,s)$ et de poids $\bm$. Si $(i,j) \prec (k,l)$, alors posons
\begin{align*}
a^c(\Lambda)_{ij,kl} &= \min(\Lambda_i^{(j)}, \Lambda_k^
{(l)}) - \Phi_i^{(j)}, \\
b^c(\Lambda)_{ij} &= \Lambda_i^{(j)} - \Phi_i^{(j)}.
\end{align*}
\'Evidemment, ces entiers d\'ependent du choix d'un pr\'esymbole repr\'esentant pour $\Lambda$. Par contre, $b^c(\Lambda)_{ij}$ est ind\'ependant de $r$ et de $s$ dans le sens suivant: si l'on change $r$ ou $s$ mais garde
la m\^eme forme pour le pr\'esymbole repr\'esentant, alors $b^c(\Lambda)_{ij}$ ne change pas.

En effet, il est clair que $b^c(\Lambda)_{ij}$ n'est autre que $\balpha_i^{(j)}$, mais cette notation sera quand m\^eme utile: elle permet de r\'ef\'erer aux coefficients d'une multipartition en termes des positions d'un symbole, au lieu de fixer une num\'erotation de ses coefficients \`a l'avance.

Enfin, il est \`a noter que l'on ne d\'efinit pas $a^c(\Lambda)_{ij,kl}$ si 
$(i,j) \not\prec (k,l)$. Ensuite, on pose
\[
a^c(\Lambda) = \sum_{(i,j) \prec (k,l)} a^c(\Lambda)_{ij,kl}
\qquad\text{et}\qquad
b^c(\Lambda) = \sum_{(i,j) \prec (k,l)} b^c(\Lambda)_{ij}
\]
Il est facile de v\'erifier que $a^c(\Lambda)$ et $b^c(\Lambda)$ sont bien d\'efinis. De plus, comme on a remarqu\'e ci-dessus, $b^c(\Lambda)$ est ind\'ependant de $r$ et de $s$. Il est clair que
\[
a^c(\Lambda)_{ij,kl} \le b^c(\Lambda)_{ij},
\qquad\text{avec \'egalit\'e si et seulement si $\Lambda_i^{(j)} \le \Lambda_k^{(l)}$}
\]
pour tout symbole $\Lambda$.  En comparant avec~\eqref{eqn:dist-defn}, on obtient imm\'ediatement le r\'esultat suivant:

\begin{lem}\label{lem:a-b-dist-combin}
Soit $\Lambda$ un symbole. On a que $a^c(\Lambda) \le b^c(\Lambda)$, avec \'egalit\'e si et seulement si $\Lambda$ est distingu\'e.
\end{lem}

Nous remarquons aussi que $b^c$ est ``additif'' (tandis que $a^c$ ne l'est pas 
en g\'en\'eral): soient $\bbeta' \in \cP(e,1,n')$ et $\bbeta'' \in \cP(e,1,n'')$,
et posons $\balpha = \bbeta' + \bbeta'' \in \cP(e,1,n'+n'')$. Fixons deux 
entiers $r$ et $s$ ainsi qu'un poids $\bm$. Il est clair que
\[
b^c(\Lambda^{r,s}_\bm(\balpha))_{ij} = b^c(\Lambda^{r,s}_\bm(\bbeta'))_{ij} + b^c(\Lambda^{r,s}_\bm(\bbeta''))_{ij},
\]
car cette \'egalit\'e \'equivaut au fait que $\balpha_i^{(j)} = \bbeta'_i{}^{(j)} + \bbeta''_i{}^{(j)}$. Le lemme suivant est donc \'evident.

\begin{lem}\label{lem:b-additif-combin}
Soient $\bbeta' \in \cP(e,1,n')$ et $\bbeta'' \in \cP(e,1,n'')$, et posons
$\balpha = \bbeta' + \bbeta'' \in \cP(e,1,n'+n'')$. Alors
\[
b^c(\Lambda^{r,s}_\bm(\balpha)) = b^c(\Lambda^{r,s}_\bm(\bbeta')) + b^c(\Lambda^{r,s}_\bm(\bbeta'')).
\]
\end{lem}

Comme on l'a d\'ej\`a remarqu\'e, les poids les plus importants sont $\bb$
et $\bd$. Les deux propositions suivantes permettent de comparer $a^c$ et
$b^c$ pour des symboles de ces deux poids.

\begin{prop}\label{prop:bc-rotation}
Si $\balpha \in \cP(e,1,n)$, alors $b^c(\Lambda^{r,s}_\bb(\balpha)) =
b^c(\rot(\Lambda^{r,0}_\bd(\balpha)))$.
\end{prop}
\begin{proof}
Posons $\Lambda = \Lambda^{r,s}_\bb(\balpha)$ et $\Lambda' = \Lambda^{r,0}_\bd(\balpha)$. Il sera commode d'\'ecrire les formules pour $b^c(\Lambda)$ et $b^c(\rot(\Lambda))$ dans une forme l\'eg\`erement diff\'erente. Posons
\[
c(\Lambda)_{ij} = |\{ \text{positions $(k,l)$ dans $\Lambda$ telles que
$(i,j) \prec (k,l)$}\}|;
\]
la d\'efinition de $c(\rot(\Lambda'))_{ij}$ est analogue. On a alors
\begin{equation}\label{eqn:bc-mult}
b^c(\Lambda) = \sum_{(i,j)} c(\Lambda)_{ij} b^c(\Lambda)_{ij}
\qquad\text{et}\qquad
b^c(\Lambda') = \sum_{(i,j)} c(\rot(\Lambda'))_{ij} b^c(\rot(\Lambda'))_{ij}.
\end{equation}

Pour comparer les symboles $\Lambda$ et $\Lambda'$, choisissons des
pr\'esymboles repr\'esentants tels que toute ligne contient $m$
coefficients, \`a l'exception de $\Lambda_0$, qui contient $m+1$
coefficients. Supposons, sans perte de g\'en\'eralit\'e, que $\Lambda_0^{(0)} = 0$. \'Evidemment, $\Lambda$ contient un coefficient de plus que
$\Lambda'$, dans la
position $(0,m)$, laquelle est la plus grande position dans l'ordre
$\prec$. Il est donc \'evident que
\[
c(\rot(\Lambda'))_{ij} = c(\Lambda')_{ij} = c(\Lambda)_{ij}-1.
\]
Autrement dit,
\[
c(\Lambda')_{ij} = c(\Lambda)_{s(i,j)}
\]
o\`u $s(i,j)$ d\'esigne la plus petite position dans $\Lambda$ (\`a l'\'egard de $\prec$) qui est plus grande que $(i,j)$. \'Evidemment, on a que
\begin{equation}\label{eqn:s-defn}
s(i,j) =
\begin{cases}
(e-1,j) & \text{si $i = 0$,} \\
(0,j+1) & \text{si $i = 1$,} \\
(i-1,j) & \text{si $1 < i < e$.}
\end{cases}
\end{equation}
D'autre part, si l'on pose $\Phi = \Phi^{r,s}_\bb$ et $\Phi' =
\Phi^{r,0}_\bd$, alors on a que
\begin{equation}\label{eqn:bd-compare}
\rot(\Lambda')_i{}^{(j)} =
\begin{cases}
\Lambda_{e-1}^{(j)}-s \\
\Lambda_0^{(j+1)}-r \\
\Lambda_{i-1}^{(j)}-s
\end{cases}
\quad\text{et}\quad
\rot(\Phi')_i^{(j)} = \Phi'_i{}^{(j)} =
\begin{cases}
\Phi_{e-1}^{(j)}-s & \text{si $i = 0$,} \\
\Phi_0^{(j+1)}-r & \text{si $i = 1$,} \\
\Phi_0^{(j)}-s & \text{si $1 < i < e$.}
\end{cases}
\end{equation}
Il est imm\'ediat que
\[
b^c(\rot(\Lambda'))_{ij} = b^c(\Lambda)_{s(i,j)}.
\]
Puisque $b^c(\Lambda)_{0,0} = 0$, on a que
\begin{multline*}
b^c(\Lambda) = \sum_{(i,j)} c(\Lambda)_{i,j} b^c(\Lambda)_{i,j}\\
= \sum_{(i,j) \ne (0,m)} c(\Lambda)_{s(i,j)} b(\Lambda)_{s(i,j)}
= \sum_{(i,j)} c(\rot(\Lambda'))_{i,j} b(\rot(\Lambda'))_{i,j}
= b^c(\rot(\Lambda')).
\end{multline*}
\end{proof}

\begin{prop}\label{prop:ac-rotation}
Si $\balpha \in \cP(e,1,n)$, alors $a^c(\Lambda^{r,r}_\bb(\balpha)) = a^c(\rot(\Lambda^{r,0}_\bd(\balpha)))$.
\end{prop}
\begin{proof}
La preuve de cette proposition est tr\`es proche de celle de la
proposition pr\'ec\'edente. Reprenons les notations de cette preuve-l\`a:
on pose $\Lambda = \Lambda^{r,r}_\bb(\balpha)$ et $\Lambda' =
\Lambda^{r,0}_\bd(\balpha)$, et on en choisit des pr\'esymboles
repr\'esentants de la mani\`ere d\'ecrite au paragraphe qui pr\'ec\`ede la
d\'efinition~\eqref{eqn:s-defn} de $s(i,j)$.

Les formules~\eqref{eqn:bd-compare}, dans le cas o\`u $s = r$, disent pr\'ecis\'ement que
\[
\rot(\Lambda')_i^{(j)} = \Lambda_{i'}^{(j')}-r
\qquad\text{et}\qquad
\rot(\Phi')_i^{(j)} = \Phi_{i'}^{(j')}-r
\]
o\`u $(i',j') = s(i,j)$. Il s'ensuit que
\[
a^c(\rot(\Lambda'))_{ij,kl} = a^c(\Lambda)_{s(i,j), s(k,l)}.
\]
On a donc
\begin{multline*}
a^c(\rot(\Lambda'))
= \sum_{(i,j) \prec (k,l)} a^c(\rot(\Lambda'))_{ij,kl} \\
= \sum_{(i,j) \prec (k,l)} a^c(\Lambda)_{s(i,j),s(k,l)}
= \sum_{\substack{(i,j) \prec (k,l) \\ (i,j) \ne (0,0)}} a^c(\Lambda)_{ij,kl}
= a^c(\Lambda),
\end{multline*}
o\`u la derni\`ere \'egalit\'e est cons\'equence du fait que $a^c(\Lambda)_{00,kl} = 0$ pour toute position $(k,l)$.
\end{proof}

Enfin, nous d\'ecrivons la relation entre les fonctions combinatoires
$a^c$ et $b^c$ et les fonctions $a$ et $b$ de la section pr\'ec\'edente.
Le r\'esultat suivant a \'et\'e \'etabli par Malle~\cite{M1}:

\begin{prop}\label{prop:fn-a}
\begin{enumerate}
\item Si $\balpha \in \cP(e,1,n)$, alors $a(\balpha) =
a^c(\Lambda^{1,0}_\bb(\balpha))$. De plus, $E_\balpha$ est sp\'eciale si et
seulement si $\Lambda^{1,0}_\bb(\balpha)$ est distingu\'e.
\item Si $\balpha \in \cP(e,e,n)$, alors $a(\balpha) =
a^c(\Lambda^{1,0}_\bd(\balpha))$. De plus, $E_{\balpha,l}$, o\`u $1 \le l
\le s_e(\balpha)$, est sp\'eciale si et seulement si
$\Lambda^{1,0}_\bd(\balpha)$ est distingu\'e. \qed
\end{enumerate}
\end{prop}

L'analogue de cet \'enonc\'e pour $b$ et $b^c$ est donn\'e ci-dessous. Il
doit \^etre possible (voire, facile) d'en donner une preuve purement
combinatoire \`a partir de la formule~\eqref{eqn:b}, mais nous
effectuerons une preuve diff\'erente plus tard. 

\begin{prop*}[Voir la Proposition~\ref{prop:Ge1n-j-ind}]
Si $\balpha \in \cP(e,1,n)$, alors 
$$b(\balpha) = b^c(\Lambda^{r,s}_\bb(\balpha))=
b^c(\rot(\Lambda^{r,0}_\bd(\balpha))).$$
\end{prop*}

Les deux corollaires suivants sont maintenant des cons\'equences
imm\'ediates des Lemmes~\ref{lem:a-b-dist-combin}
et~\ref{lem:b-additif-combin}. (Pourtant, nous \'eviterons d'utiliser
ces corollaires avant de terminer la preuve de la proposition
pr\'ec\'edente).

\begin{cor}\label{cor:a-b-dist}
Soit $\balpha \in \cP(e,1,n)$, et posons $\Lambda =
\Lambda^{r,s}_\bm(\balpha)$, o\`u $\bm$ est \'egal soit \`a $\bb$, soit \`a $\bd$.
Alors $a^c(\Lambda) \le b(\balpha)$, avec \'egalit\'e si
et seulement si $\Lambda$ est distingu\'e.\qed
\end{cor}

\begin{cor}\label{cor:b-additif}
Soient $\bbeta' \in \cP(e,1,n')$ et $\bbeta'' \in \cP(e,1,n'')$, et posons
$\balpha = \bbeta' + \bbeta'' \in \cP(e,1,n'+n'')$. Alors
$b(\balpha) = b(\bbeta') + b(\bbeta'')$.\qed
\end{cor}

%**************************************************************************
\section{L'induction tronqu\'ee}
\label{sect:tronquee}
%**************************************************************************

Soit $W$ un groupe de r\'eflexions complexes, et soit $W' \subset W$ un
sous-groupe engendr\'e par r\'eflexions. L'induction tronqu\'ee (ou
l'induction de MacDonald--Lusztig--Spaltenstein) est une op\'eration qui
associe \`a une repr\'esentation irr\'eductible de $W'$ une certaine
repr\'esentation irr\'eductible de $W$. Pourtant, cette op\'eration n'est
pas toujours d\'efinie.

\begin{defn}
Soit $W$ un groupe de r\'eflexions complexes sur l'espace vectoriel $V$.
Une repr\'esentation irr\'eductible $E$ de $W$ est dite
\emph{$j$-inductible} si elle intervient avec multiplicit\'e $1$ dans
$S^{b(E)}(V)$ (la composante de degr\'e $b(E)$ de l'alg\`ebre sym\'etrique
de $V$).
\end{defn}

\begin{prop}\label{prop:j-ind-defn}
Soit $W$ un groupe de r\'eflexions complexes, et soit $W' \subset W$ un
sous-groupe engendr\'e par r\'eflexions. Soit $E'$ une repr\'esentation
irr\'eductible $j$-inductible de $W'$, consid\'er\'ee comme sous-espace de
$S^{b(E')}(V)$, et soit $E$ le plus petit sous-espace de $S^{b(E')}(V)$
qui contient $E'$ et est $W$-stable. Alors la $W$-repr\'esen\-tation $E
\subset S^{b(E')}(V)$ est irr\'eductible. De plus, $b(E) = b(E')$, et $E$
est elle aussi $j$-inductible.
\end{prop}
\begin{proof}
La preuve de \cite[Theorem~11.2.1]{Ca} s'\'etend aux groupes de
r\'eflexions complexes.
\end{proof}

\begin{defn}
La repr\'esentation $E$ construite dans la proposition pr\'ec\'edente est
appel\'ee l'\emph{induite tronqu\'ee} de $E'$, et est not\'ee
$j_{W'}^W(E')$ (ou simplement $j(E')$ s'il n'y a aucun risque
d'ambiguit\'e).
\end{defn}

Une autre description de l'induite tronqu\'ee est comme suit: si $E'$ est
$j$-inductible, alors il y a une unique composante irr\'eductible $E$ de
l'induite (ordinaire) $\Ind_{W'}^W E'$ telle que $b(E) = b(E')$. Cette
repr\'esentation $E$ est l'induite tronqu\'ee de $E'$,
voir \cite[Proposition~11.2.5]{Ca}. (Si $E'$ n'est pas
$j$-inductible, $E$ n'est pas forc\'ement unique).

\smallskip

Les deux propositions suivantes sont bien connues.

\begin{prop}[Transitivit\'e de l'induction tronqu\'ee]
Soit $W'' \subset W' \subset W$ une suite de sous-groupes engendr\'es par
r\'eflexions, et soit $E''$ une repr\'esentation $j$-inductible de $W''$.
Alors $j_{W''}^W(E'') = j_{W'}^W(j_{W''}^{W'}(E''))$.
\end{prop}

\begin{prop}\label{prop:j-ind-tensor}
Soient $W_1$ et $W_2$ deux groupes de r\'eflexions complexes, op\'erant
sur $V_1$ et $V_2$ respectivement. Alors $W = W_1 \times W_2$ est un
groupe de r\'eflexions complexes sur $V = V_1 \oplus V_2$. Soit $E$ une
repr\'esentation irr\'eductible de $W$; elle se d\'ecompose en produit
tensoriel $E_1 \boxtimes E_2$, o\`u $E_i$ est une repr\'esentation
irr\'eductible de $W_i$, et on a que $b(E) = b(E_1) + b(E_2)$. Alors $E$ est 
$j$-inductible si et seulement si $E_1$ et $E_2$ le sont.
\end{prop}

\begin{lem}\label{lem:signe-j-ind}
Soit $\epsilon$ le caract\`ere signe du groupe sym\'etrique $\fS_n$,
consid\'er\'e comme repr\'esentation de $G(e,1,n)$. Posons $E = \epsilon
\otimes \gamma_e^k$. Alors le degr\'e fant\^ome de $E$ est
\[
R_E(X) = X^{kn + e(n^2 - n)/2}.
\]
En particulier, cette repr\'esentation est $j$-inductible.
\end{lem}
\begin{proof}
La partition de $n$ qui correspond au caract\`ere signe de $\fS_n$ est $(1
\le \cdots \le 1)$, et donc la multipartition de $\epsilon \otimes
\gamma_e^k$ est $\balpha = (\balpha_0, \ldots, \balpha_{e-1})$, o\`u
\[
\balpha_i =
\begin{cases}
(1 \le \cdots \le 1) & \text{si $i = k$,} \\
(0) & \text{si $i \ne k$.}
\end{cases}
\]
On a donc
\[
A_i(\balpha) =
\begin{cases}
\{1,2, \ldots, n\} & \text{si $i = k$,} \\
\{0\} & \text{si $i \ne k$.}
\end{cases}
\]
\'Evidemment, $\Delta(A_i(\balpha),X) = \Theta(A_i(\balpha),X) = 1$ si $i
\ne k$. De plus, $n_i(\balpha) = c_i(\balpha) = 0$ si $i \ne k$, tandis
que $n_k(\balpha) = n$ et
\begin{multline*}
c_k(\balpha) = \sum_{l = 0}^{n-1} \binom{l}{2}
= \frac{1}{2}\sum_{l=0}^{n-1}(l^2 - l) \\
= \frac{1}{2}\left(\frac{(n-1)n(2n-1)}{6} - \frac{(n-1)n}{2}\right)
= \frac{n^3 - 3n^2 + 2n}{6}.
\end{multline*}
La formule~\eqref{eqn: fantome} se r\'eduit donc \`a
\begin{equation}\label{eqn:j-ind-lem}
R_E(X) = \prod_{h=1}^n (X^{eh} - 1) \cdot
\frac{\Delta(A_k(\balpha),X^e)}{\Theta(A_k(\balpha),X^e)} \cdot
X^{kn - e(n^3 -3n^2 +2n)/6}.
\end{equation}
Ensuite, on a que
\[
\begin{array}{r@{{}={}}c@{{}={}}l}
\Delta(A_k(\balpha),X^e) &
\displaystyle\prod_{\substack{a, b \in A_k(\balpha) \\ b < a}} (X^{ea} - X^{eb})
& \displaystyle\prod_{j = 1}^n \prod_{l = j+1}^n (X^{el} - X^{ej}), \\
\Theta(A_k(\balpha),X^e) &
\displaystyle\prod_{a \in A} \prod_{l=1}^a (X^{el} - 1) 
& \displaystyle\prod_{j=1}^n \prod_{l = 1}^j
(X^{el} - 1).
\end{array}
\]
Alors, \'evidemment, on a que
\[
\frac{\Delta(A_k(\balpha),X^e)}{\Theta(A_k(\balpha),X^e)} = 
\prod_{j=1}^n B_j
\qquad\text{o\`u}\qquad
B_j = \frac{\prod_{l=j+1}^n (X^{el} - X^{ej})}{\prod_{l=1}^j (X^{el}-1)}.
\]
Si $j < n$, alors
\begin{multline*}
B_j = \frac{(X{e(j+1)}-X^{ej})(X^{e(j+2)} - X^{ej}) \cdots (X^{en} - X^{ej})}
{(X^{e} - 1)(X^{2e} - 1) \cdots (X^{ej} - 1)} \\
= \frac{X^{ej(n-j)}\cdot (X^{e} - 1)(X^{2e}-1) \cdots (X^{e(n-j)} - 1)}
{(X^{e} - 1)(X^{2e} - 1) \cdots (X^{ej} - 1)}.
\end{multline*}
Si l'on pose $c_j = X^{ej(n-j)}$ et $B'_j = c_j^{-1} B_j$, alors il est
\'evident que $B'_j B'_{n-j} = 1$ si $1 \le j < n$; de plus, dans le cas
o\`u $n$ est pair, on a que $B'_{n/2} = 1$. Par cons\'equent,
\[
\prod_{i=1}^{n-1} B_j = \prod_{i=1}^{n-1} c_j B'_j = \prod_{i=1}^{n-1} c_j.
\]
D'autre part, on a que $B_n = 1/\prod_{l = 1}^n (X^{el} - 1)$.
La formule~\eqref{eqn:j-ind-lem} s'\'ecrit donc
\begin{align*}
R_E(X) &= \prod_{h=1}^n (X^{eh} - 1) \cdot \prod_{i=1}^n B_j \cdot X^{kn -
e(n^3 - 3n^2 +2n)/6} \\
&= \prod_{h=1}^n (X^{eh} - 1) \cdot \prod_{i=1}^{n-1} c_j \cdot B_n \cdot
X^{kn - e(n^3 - 3n^2 +2n)/6} \\
&= \prod_{i=1}^{n-1} c_j \cdot X^{kn - e(n^3 - 3n^2 +2n)/6}
= X^{\sum_{i=1}^{n-1} ej(n-j) + kn - e(n^3 - 3n^2 + 2n)/6}.
\end{align*}
Il est facile de v\'erifier que
\[
\sum_{i=1}^{n-1} ej(n-j) = en \sum_{i=1}^{n-1} j - e \sum_{i=1}^{n-1} j^2
= \frac{e(n^3 -n)}{6},
\]
et donc
$R_E(X) = X^{kn + e(n^2 -n)/2}$.
\end{proof}

En particulier, on voit que $b(\epsilon \otimes \gamma_e^k) = kn + e(n^2 -
n)/2$. La l\'eg\`ere g\'en\'eralisation suivante est imm\'ediate:

\begin{cor}\label{cor:signe-b}
Soit $\bn = (n_0, \ldots, n_m)$ une suite d'entier positifs, et posons
$n = n_0 + \cdots +n_m$ et $G(e,1,\bn) = G(e,1,n_0) \times \cdots \times
G(e,1,n_m)$. Soit $E$ la
repr\'esentation $\epsilon \otimes \gamma_e^k$ de $G(e,1,\bn)$.  Alors
\[
b(E) = kn + \frac{e}{2} \sum_{i=0}^m (n_i^2 - n_i).
\]
\end{cor}

\begin{lem}\label{lem:Ge1n-j-ind}
Soit $\alpha$ une partition de $n$, et soit $\alpha^* = (0 \le
\alpha^{*(0)} \le \cdots \le \alpha^{*(m)})$ sa partition duale. Posons
$\fS_{\alpha^*} = \fS_{\alpha^{*(0)}} \times \cdots \times
\fS_{\alpha^{*(m)}}$ et $G(e,1,\alpha^*) = G(e,1,\alpha^{*(0)}) \times
\cdots \times G(e,1,\alpha^{*(m)})$. Soit $\epsilon$ le caract\`ere signe
de $\fS_{\alpha^*}$, consid\'er\'e comme caract\`ere de $G(e,1,\alpha^*)$
via la projection naturelle $G(e,1,\alpha^*) \to \fS_{\alpha^*}$. Alors
\[
j_{G(e,1,\alpha^*)}^{G(e,1,n)} (\epsilon \otimes \gamma_e^k) = E_\alpha
\otimes \gamma_e^k.
\]
\end{lem}
\begin{proof}
Rappelons les notations de la Section~\ref{subsection: repr}: le groupe $G(e,1,n) = (\Z/e\Z)^n \rtimes \fS_n$ agit sur $V = \bigoplus_{i=1}^n \C e_i$. On a d\'efini un ensemble de r\'eflexions $\{t, s_1, \ldots, s_{n-1}\}$ qui engendre $G(e,1,n)$. Les $s_i$ engendre le sous-groupe $\fS_n$. D'autre part, posons $T = (\Z/e\Z)^n$. $T$ est le plus petit sous-groupe distingu\'e contenant la r\'eflexion $t$.

Posons $\sigma_h = \sum_{l=0}^h \alpha^{*(l)}$, ainsi que $\sigma_{-1} =
0$. Consid\'erons les \'el\'ements suivants de l'alg\`ebre sym\'etrique sur
$V$:
\begin{align}
P_h &= \prod_{\sigma_{h-1} < i < j \le \sigma_h} (e_i - e_j)  &
Q_h &= \prod_{\sigma_{h-1} < i < j \le \sigma_h} 
\prod_{l=1}^{e-1} (e_i - \zeta_e^l e_j) \label{eqn:PhQh-defn} \\
P &= \prod_{h=0}^m P_h &
Q &= \prod_{h=0}^m Q_h \notag
\end{align}
ainsi que
\[
R = (e_1 \cdots e_n)^k.
\]

Il est clair que $\fS_{\alpha^*}$ agit sur $\C\cdot P$ par le caract\`ere signe, et que $Q$ et $R$ sont $\fS_{\alpha^*}$-invariants. D'autre part, le produit
\[
PQ = \prod_{\substack{0 \le h \le m \\ \sigma_h-1 < i < j \le
\sigma_h}} (e_i^e - e_j^e)
\]
est $T$-invariant, tandis que l'action de $T$ sur $\C\cdot R$ est donn\'ee
par le caract\`ere $\gamma_e^k$. En r\'esum\'e, 
\[
\begin{array}{r@{\text{ agit sur }}l@{\text{ par }}l}
\fS_{\alpha^*} & \C \cdot P & \epsilon \\
G(e,1,\alpha^*) & \C \cdot PQ & \epsilon \\
G(e,1,\alpha^*) & \C \cdot PQR & \epsilon \otimes \gamma_e^k
\end{array}
.
\]

Nous d\'emontrons maintenant que $\deg PQR = b(\epsilon \otimes
\gamma_e^k)$. Il s'ensuivra que, pour calculer l'induction tronqu\'ee de
$\epsilon \otimes \gamma_e^k$, il suffit d'\'etudier explicitement
l'action de $G(e,1,n)$ sur $PQR$. Pour commencer, consid\'erons $P_h$: son
degr\'e \'egale le nombre de valeurs distinctes que prend le couple
$(i,j)$ dans la formule~\eqref{eqn:PhQh-defn}. Puisque $\sigma_h -
\sigma_{h-1} = \alpha^{*(h)}$, on voit qu'il y a $\binom{\alpha^{*(h)}}{2}$
valeurs possibles de ce couple-l\`a. Donc
\[
\deg P_h = \frac{(\alpha^{*(h)})^2 - \alpha^{*(h)}}{2}
\qquad\text{et}\qquad
\deg Q_h = (e-1)\cdot\frac{(\alpha^{*(h)})^2 - \alpha^{*(h)}}{2}.
\]
$R$ est \'evidemment de degr\'e $kn$, et donc
\[
\deg PQR = \deg R + \sum_{h=0}^m \deg P_h Q_h
= kn + e \sum_{h=0}^m \frac{(\alpha^{*(h)})^2 - \alpha^{*(h)}}{2}.
\]
Selon le Corollaire~\ref{cor:signe-b}, il est bien le cas que $\deg PQR =
b(\epsilon \otimes \gamma_e^k)$. En sp\'ecialisant au cas o\`u $k = 0$, on
voit aussi que $\deg PQ = b(\epsilon)$, o\`u ici on calcule $b$ \`a
l'\'egard du groupe $G(e,1,\alpha^*)$.

On peut \'egalement calculer $b(\epsilon)$ pour le groupe
$\fS_{\alpha^*}$: dans ce cas, il est bien connu que $b(\epsilon) =
\sum((\alpha^{*(h)})^2 - \alpha^{*(h)})/2 = \deg P$.

Les trois polyn\^omes $P$, $PQ$, et $PQR$ ont donc les bons degr\'es pour
permettre de calculer certaines induites tronqu\'ees.  Posons
\[
\begin{gathered}
F_1 = j_{\fS_{\alpha^*}}^{\fS_n} \epsilon \subset S^{\deg P}(V) \\
F_2 = j_{G(e,1,\alpha^*)}^{G(e,1,n)} \epsilon \subset S^{\deg PQ}(V) \\
E = j_{G(e,1,\alpha^*)}^{G(e,1,n)} (\epsilon \otimes \gamma_e^k) \subset
S^{\deg PQR}(V)
\end{gathered}
\]
$F_1$ est donc le plus petit $\fS_n$-sous-module de $S^{\deg P}(V)$ qui
contient $P$, et ainsi de suite.  En particulier, $PQ \in F_2$, et $PQR
\in E$.

Il est bien connu que $F_1$ n'est autre que $E_\alpha$. Ensuite, on peut
d\'efinir une application $\fS_n$-\'equivariante $\phi: S^{\deg P}(V) \to
S^{\deg PQ}(V)$ par $\phi(f) = Qf$. Il est clair que $\phi(F_1)$ est
un sous-espace $\fS_n$-stable de $F_2$. D'autre part, puisque $PQ$ est
$T$-stable, il faut que $T$ agisse trivialement sur $F_2$, et par
cons\'equent, $F_2$ est une repr\'esentation irr\'eductible de $G(e,1,n)/T \simeq \fS_n$. On conclut que $\phi(F_1) = F_2$. Puisque la
repr\'esentation de $G(e,1,n)$ sur $F_2$ est isomorphe \`a $E_\alpha$, on
voit que sa repr\'esentation sur $E = R \cdot F_2$ \'egale $E_\alpha
\otimes \gamma_e^k$.
\end{proof}

\begin{prop}\label{prop:Ge1n-j-ind}
Toute repr\'esentation irr\'eductible de $G(e,1,n)$ est $j$-inductible. De
plus, pour toute $\balpha \in \cP(e,1,n)$, on a que $b(\balpha) =
b^c(\Lambda^{r,s}_\bb(\balpha))$.
\end{prop}
\begin{proof}
Soit $\balpha = (\balpha_0, \ldots, \balpha_{e-1}) \in \cP(e,1,n)$, et pour
chaque $i$, soit $n_i$ la somme de la partition $\balpha_i$. Posons $\bn =
(n_0, \ldots, n_{e-1})$.

Selon le Lemme~\ref{lem:Ge1n-j-ind}, chaque repr\'esentation $E_{\balpha_i}
\otimes \gamma_e^i$ est une repr\'esentation $j$-inductible de
$G(e,1,n_i)$, et donc, par la Proposition~\ref{prop:j-ind-tensor}, la
repr\'esentation
\[
F = (E_{\balpha_0} \otimes \gamma_e^0) \boxtimes \cdots \boxtimes
(E_{\balpha_{e-1}} \otimes \gamma_e^{e-1})
\]
de $G(e,1,\bn)$ est elle aussi $j$-inductible. Son induite tronqu\'ee fait
partie de son induite ordinaire, mais d'autre part, on sait que son induite
ordinaire est d\'ej\`a irr\'eductible (c'est $E_{\balpha}$). Il s'ensuit
que $E_\balpha$ est l'induite tronqu\'ee de $F$. Selon la
Proposition~\ref{prop:j-ind-defn}, la repr\'esentation $E_\balpha$ est donc 
elle aussi $j$-inductible.

Il reste \`a \'etablir la formule pour $b(\balpha)$. L'argument
pr\'e\-c\'e\-dent, com\-bi\-n\'e au Lem\-me~\ref{lem:Ge1n-j-ind}, montre
que
\[
E_\balpha = j_{G(e,1,\bn)}^{G(e,1,n)}\left( \bigboxtimes_{i=0}^{e-1} 
  j_{G(e,1,\balpha_i^*)}^{G(e,1,n_i)} (\epsilon \otimes \gamma_e^i)
\right).
\]
Au vu de la Proposition~\ref{prop:j-ind-tensor} et de l'additivit\'e de
la fonction $b^c$ (voir les commentaires qui suivent le
Lemme~\ref{lem:a-b-dist-combin}), il suffit d'\'etablir la formule pour
$b(\balpha)$ dans le cas o\`u $E_\balpha = \epsilon
\otimes \gamma_e^k$. Dans ce cas, toutes les partitions $\balpha_i$ sont
nulles si $i \ne k$, tandis que
\[
\balpha_k = (0 \le \underbrace{1 \le \cdots \le 1}_{\text{$n$ parties}}).
\]
Posons $\Lambda = \rot(\Lambda^{r,0}_\bd(\balpha))$ (o\`u $r$ est un
entier positif quelconque). Supposons, en plus, que chaque ligne de
$\Lambda$ contienne $n$ coefficients (num\'erot\'es de $0$ \`a $n-1$).
Si l'on pose
\[
k' =
\begin{cases}
k+1 & \text{si $0 \le k < e-1$,} \\
0   & \text{si $k = e-1$,}
\end{cases}
\]
alors on a que
\[
b^c(\Lambda)_{ij} =
\begin{cases}
1 & \text{si $j = k'$,} \\
0 & \text{sinon}
\end{cases}
\qquad\text{et}\qquad
c(\Lambda)_{ij} =
\begin{cases}
e(n-j-1) + (i-1) \text{si $i > 0$,} \\
e(n-j-1) + (e-1) \text{si $i = 0$.}
\end{cases}
\]
En particulier, on voit que
\[
c(\Lambda)_{k'j} = e(n-j-1) + k.
\]
La formule~\eqref{eqn:bc-mult} donne donc que
\begin{multline*}
b^c(\Lambda) = 
\sum_{(i,j)} c(\Lambda)_{ij}b^c(\Lambda)_{ij} =
\sum_{j=0}^{n-1} c(\Lambda)_{k'j} = \sum_{j=0}^{n-1}
(e(n-j-1) + k) \\
= e\left(n^2 - \sum_{j=0}^{n-1} j - n\right) + kn = en^2 -
e\frac{n(n-1)}{2} - en + kn\\
 = kn + e(n^2-n)/2.
\end{multline*}
Il d\'ecoule du Corollaire~\ref{cor:signe-b} que $b^c(\Lambda) =
b(\epsilon \otimes \gamma_e^k)$.
\end{proof}

\begin{rem} \label{rem: j}
La preuve de la proposition ci-dessus montre que toute re\-pr\'e\-sen\-ta\-tion
ir\-r\'eductible de $G(e,1,n)$ est de la forme
$$j_{G(e,1,\balpha_0^*)\times\cdots\times G(e,1,\balpha_{e-1}^*)}^{G(e,1,n)}
\left((\epsilon\otimes\gamma_e^0)\otimes\cdots\otimes
(\epsilon\otimes\gamma_e^{e-1})\right),$$
o\`u, pour chaque $i\in\{0,1,\ldots,e-1\}$, on a not\'e 
$\balpha_i^*=(0\le\balpha_i^{*(0)}\le\cdots\le\balpha_i^{*(m)})$ 
la partition duale de la partition $\balpha_i$ et o\`u
$\balpha=(\balpha_0,\ldots,\balpha_{e-1})\in\cP(e,1,n)$. Ceci constitue une
g\'en\'eralisation au groupe $G(e,1,n)$ du r\'esultat connu pour les
groupes de Weyl de type $B_n$ (voir \cite[Proposition 11.4.2]{Ca} ou
\cite{L}).
\end{rem}
 
\begin{prop}\label{prop:Geen-j-ind}
Si $\balpha \in \cP(e,e,n)$ et $\Lambda^{r,0}_\bd(\balpha)$ est distingu\'e, alors toutes les repr\'esentations $E_{\balpha,l}$ ($1 \le l \le s_e(\balpha)$) sont $j$-inductibles.
\end{prop}
\begin{proof}
Soit $\balpha \in \cP(e,e,n)$ une multipartition telle que les repr\'esentations $E_{\balpha,l}$ ($1 \le l \le s_e(\balpha)$) soient sp\'eciales. (Il est clair que la propri\'et\'e d'\^etre sp\'ecial ne d\'epend pas de $l$, puisque les fonctions $a$ et $b$ ne d\'ependent que de $\balpha$).

Si $\bbeta \in \cP(e,1,n)$, on sait que chaque $E_{\balpha,l}$ intervient dans 
la restriction de $E_\bbeta$ \`a $G(e,e,n)$ si et seulement si $\bbeta$ est une 
multipartition au-dessus de $\balpha$ (et dans ce cas, $E_{\balpha,l}$ intervient 
dans $E_\bbeta$ avec multiplicit\'e $1$). En particulier, il s'ensuit que
\[
b(\balpha) = \min \{ b(\bbeta) \mid
\text{$\bbeta \in \cP(e,1,n)$ est au-dessus de $\balpha$} \}.
\]
Soit $\tilde\balpha \in \cP(e,1,n)$ une multipartition au-dessus de $\balpha$ 
telle que $b(\balpha) = b(\tilde\balpha)$. On voit que $E_{\balpha,l}$ est 
$j$-inductible si et seulement si $\tilde\balpha$ est l'{\it unique} multipartition 
au-dessus de $\balpha$ en laquelle la valeur de la fonction $b$ \'egale 
$b(\balpha)$. 

Il est clair que les autres multipartitions au-dessus de $\balpha$ s'obtiennent \`a partir de $\tilde\balpha$ par rotation. Donc nous voudrions d\'emontrer que
\[
b(\rot^k(\tilde\balpha)) > b(\balpha)
\qquad\text{si $\rot^k(\tilde\balpha) \ne \tilde\balpha$}.
\]

Posons $\Lambda = \Lambda^{r,0}_\bd(\balpha)$ et $\tilde\Lambda = \Lambda^{r,0}_\bd(\tilde\balpha)$. Il est clair que
\[
a^c(\Lambda) = a^c(\tilde\Lambda) = a^c(\rot(\tilde\Lambda)).
\]
D'une part, puisque $\Lambda$ est distingu\'e, on sait que $b(\balpha) =
a(\balpha) = a^c(\Lambda)$. D'autre part, la
Proposition~\ref{prop:Ge1n-j-ind} dit que $b(\tilde\balpha) =
b^c(\rot(\tilde\Lambda))$. On conclut que $a^c(\rot(\tilde\Lambda)) =
b^c(\rot(\tilde\Lambda))$, et donc, selon le
Lemme~\ref{lem:a-b-dist-combin}, que $\rot(\tilde\Lambda)$ est distingu\'e.

S'il y avait une multipartition $\rot^k(\tilde\balpha)$, diff\'erente de
$\tilde\balpha$, telle que $b(\rot^k(\tilde\balpha)) = b(\balpha)$, alors
on saurait que $b^c(\rot^{k+1}(\tilde\Lambda)) = b(\balpha)$ aussi.
Pourtant, on sait que $a^c(\rot^{k+1}(\tilde\Lambda)) =
a^c(\rot(\tilde\Lambda))$, et donc on voit que $\rot^{k+1}(\tilde\Lambda)$
devrait \^etre distingu\'e. De plus, puisque $\tilde\balpha \ne
\rot^k(\tilde\balpha)$, on sait que $\rot(\tilde\Lambda) \ne
\rot^{k+1}(\tilde\Lambda)$. Mais il est \'evident qu'une $\rot$-orbite de
symboles de poids $\bd$ poss\`ede au plus un membre distingu\'e.

Ainsi, $\tilde\balpha$ est bien l'unique multipartition au-dessus de $\balpha$ telle que $b(\tilde\balpha) = b(\balpha)$, et donc $E_{\balpha,l}$ est $j$-inductible.
\end{proof}

Au cours de la preuve de la proposition pr\'ec\'edente, nous avons \'etabli
le fait suivant: si $\balpha \in \cP(e,e,n)$ est une multipartition telle
que $\Lambda^{r,0}_\bb(\balpha)$ est distingu\'e, et si l'on d\'efinit
$\tilde\balpha \in \cP(e,1,n)$ par l'\'equation $E_{\tilde\balpha} =
j_{G(e,e,n)}^{G(e,1,n)} E_{\balpha,l}$, alors
$\rot(\Lambda^{r,0}_\bd(\tilde\balpha))$ est distingu\'e. Les
Propositions~\ref{prop:bc-rotation} et~\ref{prop:ac-rotation} impliquent
alors que $a^c(\Lambda^{r,r}_\bb(\tilde\balpha)) =
b^c(\Lambda^{r,r}_\bb(\tilde\balpha))$, et donc que
$\Lambda^{r,r}_\bb(\tilde\balpha)$ est \'egalement distingu\'e. Cette
observation fait partie du corollaire suivant.

\begin{cor}\label{cor:Geen-Ge1n}
Soient $\balpha \in \cP(e,e,n)$ et $\tilde\balpha \in \cP(e,1,n)$. Les trois conditions suivantes sont \'equivalentes:
\begin{enumerate}
\item $\Lambda^{r,0}_\bd(\balpha)$ est distingu\'e, et $E_{\tilde\balpha} \simeq j_{G(e,e,n)}^{G(e,1,n)} E_{\balpha,l}$ pour tout $l$, $1 \le l \le s_e(\balpha)$.\label{it:1}
\item $\Lambda^{r,r}_\bb(\tilde\balpha)$ est distingu\'e, et $\tilde\balpha$ est au-dessus de $\balpha$.\label{it:2}
\item $\rot(\Lambda^{r,0}_\bd(\tilde\balpha))$ est distingu\'e, et $\tilde\balpha$ est au-dessus de $\balpha$.\label{it:3}
\end{enumerate}
\end{cor}
\begin{proof}
Il reste \`a montrer que les conditions~\eqref{it:2} et~\eqref{it:3} sont
\'equivalentes et qu'elles impliquent la condition~\eqref{it:1}. Soit
$\tilde\balpha \in \cP(e,1,n)$, et soit $\balpha$ son image dans
$\cP(e,1,n)$. Il d\'ecoule des Propositions~\ref{prop:bc-rotation}
et~\ref{prop:ac-rotation} que $\Lambda^{r,r}_\bb(\tilde\balpha)$ est
distingu\'e si et seulement si $\rot(\Lambda^{r,0}_\bd(\tilde\balpha))$
l'est. Supposons que ces deux conditions soient satisfaites. On sait, par
d\'efinition, que $\Lambda^{r,0}_\bd(\balpha)$ est distingu\'e (car
$\rot(\Lambda^{r,0}_\bd(\tilde\balpha))$ l'est). Par suite, la proposition
pr\'ec\'edente nous dit que tous les $E_{\balpha,l}$ sont $j$-inductible.
En effet, leur induite tronqu\'ee (commune) doit \^etre
$E_{\tilde\balpha}$: si l'on d\'efinit $\tilde\balpha'$ par
$E_{\tilde\balpha'} = j_{G(e,e,n)}^{G(e,1,n)} E_{\balpha,l}$, alors on sait
que $\tilde\balpha'$ est au-dessus de $\balpha$ et donc est une rotation de
$\tilde\balpha$; mais on sait aussi que
$\rot(\Lambda^{r,0}_\bd(\tilde\balpha'))$ est distingu\'e. Ce dernier
\'etant une rotation du symbole distingu\'e
$\rot(\Lambda^{r,0}_\bd(\tilde\balpha))$, on voit que les deux doivent
\^etre \'egaux, et donc que $\tilde\balpha' = \tilde\balpha$.
\end{proof}

%**************************************************************************
\section{Symboles sp\'etsiaux et induction tronqu\'ee pour $G(e,1,n)$}
\label{sect:Ge1n}
%**************************************************************************

Les r\'esultats principaux de cette section (les Th\'eor\`emes~\ref{thm:Ge1n-springer} et~\ref{thm:Geen-springer}) fournissent un lien entre les symboles distingu\'es, l'induction tronqu\'ee, et les repr\'esentations sp\'eciales. Nous aurons besoin de la proposition suivante, utile pour le calcul des induites tronqu\'ees.

\begin{prop}\label{prop:Ge1n-j-ind-formule}
Soient $\bbeta' \in \cP(e,1,n')$ et $\bbeta'' \in \cP(e,1,n'')$, et posons
$\balpha = \bbeta' + \bbeta'' \in \cP(e,1,n'+n'')$. Alors
\[
j_{G(e,1,n') \times G(e,1,n'')}^{G(e,1,n'+n'')}(E_{\bbeta'} \boxtimes
E_{\bbeta''}) = E_\balpha.
\]
\end{prop}
\begin{proof}
Si $\bbeta' = (\bbeta'_0, \ldots, \bbeta'_{e-1})$ et $\bbeta'' = (\bbeta''_0, \ldots, \bbeta''_{e-1})$, soit $n'_i$ (resp.~$n''_i$) la somme de la partition $\bbeta'_i$ (resp.~$\bbeta''_i$), o\`u $0 \le i < e$. Posons aussi $\bn' = (n'_0, \ldots, n'_{e-1})$ et $\bn'' = (n''_0, \ldots, n''_{e-1})$. Comme on a remarqu\'e au cours de la preuve de la Proposition~\ref{prop:Ge1n-j-ind}, on a que
\[
E_{\bbeta'} = j_{G(e,1,\bn')}^{G(e,1,n')} (E_{\bbeta'_0} \boxtimes (E_{\bbeta'_1} \otimes \gamma_e) \boxtimes \cdots \boxtimes (E_{\bbeta'_{e-1}} \otimes \gamma_e^{e-1})),
\]
et de m\^eme pour $E_{\bbeta''}$.

Posons $n = n' + n''$, $n_i = n'_i + n''_i$, et $\bn = \bn' + \bn''$. En utilisant les d\'efinitions de $E_{\bbeta'}$ et de $E_{\bbeta''}$ et la transitivit\'e de l'induction tronqu\'ee, on trouve que
\begin{align}
&j_{G(e,1,n') \times G(e,1,n'')}^{G(e,1,n)} (E_{\bbeta'} \boxtimes
E_{\bbeta''}) \notag \\
&\qquad = j_{G(e,1,n') \times G(e,1,n'')}^{G(e,1,n)}
j_{G(e,1,\bn') \times G(e,1,\bn'')}^{G(e,1,n') \times G(e,1,n'')}
\bigboxtimes_{i=0}^{e-1}(E_{\bbeta'_i} \otimes \gamma_e^i) \boxtimes
\bigboxtimes_{i=0}^{e-1}(E_{\bbeta''_i} \otimes \gamma_e^i) \notag \\
&\qquad= j_{G(e,1,\bn)}^{G(e,1,n)} 
\bigboxtimes_{i=0}^{e-1} j_{G(e,1,n'_i) \times G(e,1,n''_i)}^{G(e,1,n_i)}
(E_{\bbeta'_i} \otimes \gamma_e^i) \boxtimes
(E_{\bbeta''_i} \otimes \gamma_e^i). \label{eqn:j-ind-somme-mp}
\end{align}
\'Etudions maintenant les facteurs du grand produit tensoriel ci-dessus: selon le Lemme~\ref{lem:Ge1n-j-ind}, pour chaque $i$, on a
\[
j_{G(e,1,n'_i) \times G(e,1,n''_i)}^{G(e,1,n_i)} (E_{\bbeta'_i} \otimes \gamma_e^i) \boxtimes (E_{\bbeta''_i} \otimes \gamma_e^i)
= j_{G(e,1,(\bbeta'_i)^*) \times G(e,1,(\bbeta''_i)^*}^{G(e,1,n_i)}
(\epsilon \otimes \gamma_e^i) \boxtimes (\epsilon \otimes \gamma_e^i)
\]
Maintenant, on remarque que la r\'eunion des parties de $(\bbeta'_i)^*$ et 
de $(\bbeta''_i)^*$ n'est autre que l'ensemble de parties de
$(\bbeta'_i + \bbeta''_i)^*$. En particulier, on a que
\[
G(e,1,(\bbeta'_i)^*) \times G(e,1,(\bbeta''_i)^*) \simeq G(e,1,(\bbeta_i)^*).
\]
Il est clair que, sous cet isomorphisme, la repr\'esentation $(\epsilon \otimes \gamma_e^i) \boxtimes (\epsilon \otimes \gamma_e^i)$ s'identifie avec la repr\'esentation $\epsilon \otimes \gamma_e^i$ de $G(e,1,(\bbeta_i)^*)$.

La formule~\eqref{eqn:j-ind-somme-mp} devient donc:
\begin{multline*}
j_{G(e,1,n') \times G(e,1,n'')}^{G(e,1,n)} (E_{\bbeta'} \boxtimes E_{\bbeta''}) 
= j_{G(e,1,\bn)}^{G(e,1,n)} \bigboxtimes_{i=0}^{e-1} j_{G(e,1,(\bbeta_i)^*)}^{G(e,1,n_i)} (\epsilon \otimes \gamma_e^i) \\
= j_{G(e,1,\bn)}^{G(e,1,n)} \bigboxtimes_{i=0}^{e-1} (E_{\bbeta_i} \otimes \gamma_e^i) = E_{\bbeta}.
\end{multline*}
\end{proof}

Ensuite, la proposition suivante, qui permet de d\'ecomposer un symbole en
somme de deux symboles plus petits, sera indispensable.

\begin{prop}\label{prop:Ge1n-somme-symbole}
Soit $\balpha \in \cP(e,1,n)$, et posons $\Lambda =
\Lambda^{r,s}_\bb(\balpha)$.  $\Lambda$ est distingu\'e si et seulement
si, pour tous les entiers positifs $r', r'', s', s''$ tels que
\begin{equation}\label{eqn:somme}
\begin{gathered}
0 \le s' \le r' \\
0 \le s'' \le r''
\end{gathered}
\qquad\text{et}\qquad
\begin{aligned}
r' + r'' &= r \\
s' + s'' &= s, 
\end{aligned}
\end{equation}
il existe des entiers $n'$, $n''$ tels que $n' + n'' = n$, et des multipartitions $\bbeta' \in \cP(e,1,n')$ et $\bbeta'' \in
\cP(e,1,n'')$ telles que
\[
\balpha = \bbeta' + \bbeta''
\qquad\text{et}\qquad
\text{$\Lambda^{r',s'}_\bb(\bbeta')$ et $\Lambda^{r'',s''}_\bb(\bbeta'')$
sont distingu\'es.}
\]
\end{prop}
\begin{proof}
Supposons d'abord qu'on a des entiers $n'$ et $n''$ et des multipartitions
$\bbeta' \in \cP(e,1,n')$ et $\bbeta'' \in \cP(e,1,n'')$ avec les
propri\'et\'es d\'ecrites ci-dessus. Posons $\Lambda' =
\Lambda^{r',s'}_\bb(\bbeta')$ et $\Lambda'' =
\Lambda^{r'',s''}_\bb(\bbeta'')$. Choisissons des pr\'esymboles
repr\'esent\-ants pour $\Lambda$, $\Lambda'$, et $\Lambda''$ de
mani\`ere que tous les trois aient la $0$-\`eme ligne \`a $m+1$
coefficients, et toutes les autres lignes \`a $m$ coefficients. Posons
aussi $\Phi = \Phi^{r,s}(\bb)$, $\Phi' = \Phi^{r',s'}(\bm)$, et $\Phi'' =
\Phi^{r'',s''}(\bm)$. D'apr\`es la d\'efinition des protosymboles, il est
\'evident que $\Phi_i^{(j)} = \Phi'_i{}^{(j)} + \Phi''_i{}^{(j)}$. Puisque
$\balpha_i^{(j)} = \bbeta'_i{}^{(j)} + \bbeta''_i{}^{(j)}$, on voit que
\[
\Lambda_i^{(j)} = \Lambda'_i{}^{(j)} + \Lambda''_i{}^{(j)}.
\]
Il en d\'ecoule que $a^c(\Lambda)_{ij,kl} = a^c(\Lambda')_{ij,kl} + a^c(\Lambda'')_{ij,kl}$, et donc que
\[
a^c(\Lambda) = a^c(\Lambda') + a^c(\Lambda'').
\]
D'autre part, on sait, d'apr\`es la proposition pr\'ec\'edente,
que
\[
E_\balpha = j_{G(e,1,n') \times G(e,1,n'')}^{G(e,1,n)} E_{\bbeta'}
\boxtimes E_{\bbeta''},
\]
et donc que $b(\balpha) = b(E_{\bbeta'} \boxtimes
E_{\bbeta''}) = b(\bbeta') + b(\bbeta'')$. Selon la
Proposition~\ref{prop:Ge1n-j-ind}, on peut conclure que
\[
b^c(\Lambda) = b^c(\Lambda') + b^c(\Lambda'').
\]
Puisque $\Lambda'$ et $\Lambda''$ sont distingu\'es, on voit que $a^c(\Lambda) = b^c(\Lambda)$ (voir le Lemme~\ref{lem:a-b-dist-combin}), et donc on d\'eduit que $\Lambda$ est distingu\'e aussi.

D'autre part, supposons maintenant que $\Lambda$ est distingu\'e. Nous voulons trouver $n'$, $n''$, $\bbeta'$, et $\bbeta''$ avec les propri\'et\'es \'enonc\'ees dans la proposition. En fait, nous allons d'abord d\'ecrire leurs symboles $\Lambda' = \Lambda^{r',s'}_\bb(\bbeta')$ et $\Lambda'' = \Lambda^{r'',s''}_\bb(\bbeta'')$, et puis nous montrerons que ces symboles-l\`a proviennent en effet des multipartitions telles que cherch\'ees.

Les deux fonctions suivantes nous seront utiles:
\begin{align*}
&\kappa': \{0, \ldots, r-1\} \to \{0 \ldots r'\} &
&\kappa'': \{0,\ldots, r-1\} \to \{0,\ldots, r''-1\} \\
&\kappa'(i) =
\begin{cases}
i  & \text{si $0 \le i \le s'$,} \\
s' & \text{si $s' \le i \le s$,} \\
i-s'' & \text{si $s \le i \le s''+r'$,} \\
r' & \text{si $s''+r' \le i \le r-1$}
\end{cases}
&
&\kappa''(i) =
\begin{cases}
0  & \text{si $0 \le i \le s'$,} \\
i-s' & \text{si $s' \le i \le s$,} \\
s'' & \text{si $s \le i \le s''+r'$,} \\
i-r' & \text{si $s''+r' \le i \le r-1$}
\end{cases}
\end{align*}
Il est clair que les fonctions $\kappa'$ et $\kappa''$ sont toutes deux
croissantes, et que
\[
\kappa'(i) + \kappa''(i) = i
\]
pour tout $i$.

Pour chaque position $(i,j)$, le coefficient $\Lambda_i^{(j)}$ s'\'ecrit
\[
\Lambda_i^{(j)} = ar + b
\qquad\text{o\`u $0 \le b \le r-1$}
\]
de mani\`ere unique. D\'efinissons $\Lambda'$ et $\Lambda''$ en posant
\[
(\Lambda')_i^{(j)} = ar' + \kappa'(b)
\qquad\text{et}\qquad
(\Lambda'')_i^{(j)} = ar'' + \kappa''(b)
\]
et puis posons $\bbeta' = \Lambda' - \Phi^{r',s'}(\bb)$ et $\bbeta'' = \Lambda'' - \Phi^{r'',s''}(\bb)$. Il faut v\'erifier que $\bbeta'$ et $\bbeta''$ sont bien des multipartitions, et que $\Lambda'$ et $\Lambda''$ sont distingu\'es.

Pour d\'emontrer que $\bbeta'$ est une multipartition, il suffit de montrer que, pour tout $i$, $\bbeta'_i{}^{(j-1)} \le \bbeta'_i{}^{(j)}$ si $j > 0$, et que $\bbeta'_i{}^{(0)} \ge 0$. Ces in\'egalit\'es \'equivalent aux suivantes:
\[
(\Lambda')_i{}^{(j-1)} \le (\Lambda')_i{}^{(j)} - r'
\qquad\text{et}\qquad
(\Lambda')_0{}^{(0)} \ge 0, \quad
\text{$(\Lambda')_i{}^{(0)} \ge s'$ si $i > 0$.}
\]
Il est \'evident que $(\Lambda')_0{}^{(0)} \ge 0$, et il est tr\`es facile d'\'etablir les autres in\'egalit\'es \`a partir des faits correspondants pour $\Lambda$:
\[
\Lambda_i^{(j-1)} \le \Lambda_i^{(j)}
\qquad\text{et}\qquad
\Lambda_0^{(0)} \ge 0, \quad 
\text{$\Lambda_i^{(0)} \ge s$ si $i > 0$.}
\]
Si l'on \'ecrit $\Lambda_i^{(j-1)} = a_1r + b_1$ et $\Lambda_i^{(j)} = a_2r
+ b_2$, alors il faut que $a_2 \ge a_1 + 1$, et, dans le cas o\`u $a_2 =
a_1+1$, on sait que $b_2 \ge b_1$. Dans le cas o\`u $a_2 > a_1 + 1$, il est
clair que $(\Lambda')_i{}^{(j-1)} \le (\Lambda')_i{}^{(j)} - r'$; par
contre, si $a_2 = a_1 + 1$, l'in\'egalit\'e cherch\'ee est cons\'equence du
fait que $\kappa'(b_2) \ge \kappa'(b_1)$. Ensuite, \'ecrivons
$\Lambda_i^{(0)} = ar + b$. Si $a > 0$, on voit que $(\Lambda')_i{}^{(0)} =
ar' + \kappa'(b) \ge r' \ge s'$; par contre, si $a = 0$, alors on sait que
$b \ge s$, et donc que $\kappa'(b) \ge s'$. Donc $\bbeta'$ est bien une
multipartition.

D\'emontrons maintenant que $\Lambda'$ est distingu\'e. Soient $(i,j)$ et
$(k,l)$ deux positions telles que $(i,j) \prec (k,l)$. \'Ecrivons
\[
\Lambda_i^{(j)} = a_1r + b_1
\qquad\text{et}\qquad
\Lambda_k^{(l)} = a_2r + b_2.
\]
On sait que $a_1r + b_1 \le a_2r + b_2$, ce qui implique que soit $a_1 <
a_2$, soit $a_1 = a_2$ et $b_1 \le b_2$. Dans tous les deux cas, on voit
que $a_1r' + \kappa'(b_1) \le a_2r' + \kappa'(b_2)$. Autrement dit,
$(\Lambda')_i^{(j)} \le (\Lambda')_k^{(l)}$, et $\Lambda'$ est distingu\'e.

Les preuves des faits que $\bbeta''$ est une multipartition et que
$\Lambda''$ est distingu\'e sont analogues.
\end{proof}

Nous pouvons maintenant \'etablir les th\'eor\`emes suivants:

\begin{thm}\label{thm:Ge1n-springer}
Soit $W = G(e,1,n)$, et soient $r$ et $s$ deux entiers tels que $0 \le s
\le r$. Posons
\[
\cS = \{ E = j_{W'}^W(E') \in \Irr(W) \mid
\text{$E'$ est une repr\'esentation sp\'eciale de $W'$} \}
\]
o\`u $W'$ parcourt les sous-groupes de $W$ de la forme
\[
\underbrace{G(e,e,n_1) \times \cdots \times G(e,e,n_s)}_{\text{$s$
facteurs}} 
\times
\underbrace{G(e,1,n_{s+1}) \times \cdots \times G(e,1,n_r)}_{\text{$r-s$
facteurs}}
\]
avec $n_1 + \cdots + n_r = n$. Alors
\[
\cS = \{ E_\balpha \mid \text{$\Lambda^{r,s}_\bb(\balpha)$ est distingu\'e}
\}.
\]
\end{thm}
\begin{proof}
Selon le Corollaire~\ref{cor:Geen-Ge1n}, pour toute multipartition $\balpha
\in \cP(e,1,n)$, $\Lambda^{1,1}_\bb(\balpha)$ est distingu\'e si et
seulement si $E_\balpha$ est l'induite tronqu\'ee d'une repr\'esenta\-tion
sp\'eciale de $G(e,e,n)$. D'autre part, on sait que $E_\balpha$ lui-m\^eme
est sp\'ecial si et seulement si $\Lambda^{1,0}_\bb(\balpha)$ est
distingu\'e. Donc $\cS$ s'\'ecrit
\[
\cS = \{ j_{W'}^W (E_{\balpha_1} \boxtimes \cdots \boxtimes E_{\balpha_r})
\},
\]
o\`u $W'$ parcourt les sous-groupes de la forme
\[
G(e,1,n_1) \times \cdots \times G(e,1,n_r),
\qquad\text{avec $n_1 + \cdots + n_r = n$,}
\]
et les multipartitions $\balpha_1, \ldots, \balpha_r$ sont telles que les
symboles
\[
\Lambda^{1,1}_\bb(\balpha_1), \ldots, \Lambda^{1,1}_\bb(\balpha_s);
\Lambda^{1,0}_\bb(\balpha_{s+1}), \ldots, \Lambda^{1,0}_\bb(\balpha_r)
\]
sont tous distingu\'es. Il s'ensuit de la
Proposition~\ref{prop:Ge1n-j-ind-formule} que $\cS$ peut \'egalement
s'\'ecrire
\[
\cS = \left\{ E_\balpha \Biggm| 
\begin{array}{c}
\text{$\balpha = \balpha_1 + \cdots + \balpha_r$, o\`u
$\Lambda^{1,1}_\bb(\balpha_1), \ldots, \Lambda^{1,1}_\bb(\balpha_s)$,} \\
\text{et $\Lambda^{1,0}_\bb(\balpha_{s+1}), \ldots,
\Lambda^{1,0}_\bb(\balpha_r)$ sont distingu\'es}
\end{array}
\right\}.
\]
Enfin, un argument par r\'ecurrence utilisant la
Proposition~\ref{prop:Ge1n-somme-symbole} montre que $E_\balpha \in \cS$ si
et seulement si $\Lambda^{r,s}_\bb(\balpha)$ est distingu\'e.
\end{proof}

\begin{thm}\label{thm:Geen-springer}
Soit $W = G(e,e,n)$, et soit $r$ un entier positif. Posons
\[
\cS = \{ E = j_{W'}^W(E') \in \Irr(W) \mid
\text{$E'$ est une repr\'esentation sp\'eciale de $W'$} \}
\]
o\`u $W'$ parcourt les sous-groupes de $W$ de la forme
\[
G(e,e,n_1) \times \cdots \times G(e,e,n_r)
\]
avec $n_1 + \cdots + n_r = n$. Alors
\[
\cS = \{ E_{\balpha,l} \mid \text{$\Lambda^{r,0}_\bd(\balpha)$ est
distingu\'e}
\}.
\]
\end{thm}
\begin{proof}
Posons $\tilde \cS = \{ j_{G(e,e,n)}^{G(e,1,n)} E_{\balpha,l} \mid
E_{\balpha,l}\in \cS \}$. Selon le Th\'eor\`eme~\ref{thm:Ge1n-springer},
$E_{\tilde\balpha} \in \tilde \cS$ si et seulement si
$\Lambda^{r,r}_\bb(\tilde\balpha)$ est distingu\'e. Ensuite, le
Corollaire~\ref{cor:Geen-Ge1n} nous dit que
$\Lambda^{r,r}_\bb(\tilde\balpha)$ est distingu\'e si et seulement si
$E_{\tilde\balpha}$ est l'induite tronqu\'ee d'un $E_{\balpha,l} \in
\Irr(G(e,e,n))$ avec $\Lambda^{1,0}_\bd(\balpha)$ distingu\'e.
\end{proof}

%**************************************************************************
\section{Sous-groupes paraboliques}
\label{sect:paraboliques}
%**************************************************************************

Le but de cette section est d'\'etablir le r\'esultat suivant:

\begin{prop}\label{prop:paraboliques}
Soit $W$ un groupe de r\'eflexions complexes sp\'etsial imprimitif, et soit 
$W' \subset W$ un sous-groupe parabolique. Si $E \in \Irr(W')$ est sp\'eciale, 
alors $j_{W'}^W E$ l'est aussi.
\end{prop}

Nous faisons d'abord le calcul de certaines induites tronqu\'ees en termes
de symboles.

\begin{prop}\label{prop:Gef1n-j-ind}
Soit $\balpha = (\balpha_0, \ldots, \balpha_{e-1})\in \cP(e,1,n)$, o\`u
\[
\balpha_i = (\balpha_i^{(0)} \le \cdots \le \balpha_i^{m_i}).
\]
Soit $f$ un entier strictement positif. L'induite 
tronqu\'ee $j_{G(e,1,n)}^{G(ef,1,n)} E_\balpha$ est isomorphe \`a
$E_\bbeta$, o\`u $\bbeta = (\bbeta_0, \ldots, \bbeta_{ef-1})$ est donn\'e
par
\[
\bbeta_{ke+i} = ( 0 \le \cdots \le \balpha_i^{(m_i -k - 2f)}
\le \balpha_i^{(m_i-k-f)} \le \balpha_i^{(m_i - k)}).
\]
\end{prop}

\begin{exm}
Il est plus facile de comprendre la proposition pr\'ec\'edente en termes
des symboles de type $(0,0)$ et l'ordre $\prec$: pour d\'eterminer
$\bbeta_i^{(j)}$, on cherche l'unique position $(k,l)$ dans $\balpha$
telle que $c_{ij}(\Lambda^{0,0}_\bb(\bbeta)) =
c_{kl}(\Lambda^{0,0}_\bb(\balpha))$, et on pose $\bbeta_i^{(j)} =
\balpha_k^{(l)}$. Par exemple, si
\[
\balpha =
\hbox{\tiny$\left(\begin{array}{@{}c@{}c@{}c@{}c@{}c@{}c@{}c@{}}
1&&2&&2&&3\\ &1&&3&&4\\ &0&&5&&6
\end{array}\right)$},
\]
alors $j_{G(3,1,27)}^{G(6,1,27)} E_\balpha \simeq E_\bbeta$, o\`u
\[
\bbeta =
\hbox{\tiny$\left(\begin{array}{@{}c@{}c@{}c@{}c@{}c@{}}
0&&2&&3\\ &1&&4\\ &0&&6\\ &1&&2\\ &0&&3\\ &0&&5
\end{array}\right)$}.
\]
\end{exm}

Pour prouver la Proposition~\ref{prop:Gef1n-j-ind}, nous aurons besoin du
lemme suivant, qui traite un cas particulier.

\begin{lem}
Supposons que $f \ge n$. La repr\'esentation $j_{G(e,1,n)}^{G(ef,1,n)}
(\epsilon \otimes \gamma_e^k)$ est isomorphe \`a $E_\bbeta$, o\`u
\[
\bbeta_j =
\begin{cases}
(1) & \text{si $j \equiv k \pmod e$ et $j < ne$,} \\
\varnothing & \text{sinon}
\end{cases}
\]
\end{lem}
Il est \`a noter que ce lemme est bien en accord avec la proposition
ci-dessus: la multipartition de la repr\'esentation $\epsilon \otimes
\gamma_e^k$ est $\balpha = (\balpha_0, \ldots, \balpha_{e-1})$, o\`u
\[
\balpha_j =
\begin{cases}
(1 \le \cdots \le 1) & \text{si $j = k$,} \\
\varnothing & \text{sinon.}
\end{cases}
\]
\begin{proof}
Posons
\begin{align*}
P &= \prod_{1 \le i < j \le n} (e_i^e - e_j^e), &
Q &= \prod_{i = 2}^n e_i^{(i-1)e}, &
R &= (e_1\cdots e_n)^k.
\end{align*}
Comme on l'a d\'emontr\'e au cours de la preuve du
Lemme~\ref{lem:Ge1n-j-ind}, $G(e,1,n)$ agit sur $\C \cdot PR$ par la
repr\'esentation $\epsilon \otimes \gamma_e^k$ de $G(e,1,n)$, et le degr\'e
de $PR$ est juste pour le calcul de son induite tronqu\'ee.

D'autre part, $E_\bbeta$ s'obtient par induction tronqu\'ee comme suit:
\begin{multline*}
E_\bbeta = j_{G(ef,1,1) \times \cdots \times G(ef,1,1)}^{G(ef,1,n)}
\left(\bigboxtimes_{h = 1}^n (\epsilon \otimes \gamma_{ef}^{k+he})\right)
\\
=  j_{G(ef,1,1) \times \cdots \times G(ef,1,1)}^{G(ef,1,n)}
\left(\bigboxtimes_{h = 0}^{n-1} \gamma_{ef}^{k+he}\right).
\end{multline*}
(On peut supprimer les $\epsilon$ figurant dans le produit tensoriel car
ils d\'esignent la repr\'esentation signe du groupe trivial $\fS_1$). On
renvoie le lecteur au Lemme~\ref{lem:Ge1n-j-ind} encore une fois pour v\'erifier que
la repr\'esentation $\bigboxtimes_{h=0}^{n-1} \gamma_{ef}^{k+he}$ est
r\'ealis\'ee par le polyn\^ome
\[
e_1^k e_2^{k+e} \cdots e_n^{k+(n-1)e} = e_2^e e_3^{2e} \cdots e_n^{(n-1)e}
\cdot (e_1 \cdots e_n)^k = QR.
\]

Pour d\'emontrer que $j_{G(e,1,n)}^{G(ef,1,n)} \simeq E_\bbeta$, il suffit
de d\'emontrer que le polyn\^ome $PR$ appartient au $G(ef,1,n)$-module
engendr\'e par $QR$.  $P$ est le produit de $n(n-1)/2$ facteurs, et
chaque terme $e_i^e$ figure dans $n-1$ d'entre eux.  Il est donc clair que
$PR$ est de la forme
\[
PR = R\sum_{\substack{0 \le h_1, \ldots, h_n \le n-1 \\ h_1+\cdots+h_n =
n(n-1)/2}} C_{h_1,\ldots,h_n} e_1^{h_1e} e_2^{h_2e} \cdots e_n^{h_ne},
\]
o\`u les $C_{h_1,\ldots,h_n} \in \Z$. Soit $w \in \fS_n \subset G(e,1,n)$
la permutation qui \'echange $i$ et $j$ et fixe tous les autres entiers. 
Alors l'action de $w$ sur $PR$ fixe les termes o\`u $h_i = h_j$ et permute
les autres termes. Mais on sait que $G(e,1,n)$ agit sur $\C \cdot PR$ par
$\epsilon \otimes \gamma_e^k$, et donc $w \cdot PR = -PR$. On en d\'eduit
que $C_{h_1,\ldots,h_n} = 0$ si $h_i = h_j$.

Autrement dit, pour que $C_{h_1, \ldots, h_n}$ soit non nul, il faut que
les $h_i$ soient des entiers positifs distincts et inf\'erieurs ou \'egaux
\`a $n-1$. \'Evidemment, \`a permutation pr\`es, la
seule possibilit\'e est $(h_1,\ldots,h_n) = (0,1,\ldots, n-1)$. En
particulier, $Q$ est l'un des termes figurant dans la somme ci-dessus, et
les autres s'en obtiennent par l'action de $\fS_n \subset G(ef,1,n)$. $PR$
et donc bien dans le $G(ef,1,n)$-module engendr\'e par $QR$.
\end{proof}

Le corollaire suivant est maintenant imm\'ediat:

\begin{cor}\label{cor:Gef1n-j-ind}
Supposons que $f \ge n$, et soit $\balpha = (\balpha_0, \ldots, \balpha_{e-1}) \in \cP(e,1,n)$, o\`u $\balpha_i = (0 \le \balpha_i^{(0)} \le \cdots \le \balpha_i^{(m_i)})$. La repr\'esentation $j_{G(e,1,n)}^{G(ef,1,n)} E_\balpha$ est isomorphe \`a $E_\bbeta$, o\`u
\[
\bbeta_j = 
\begin{cases}
(\balpha_i^{(m_i - k)}) & \text{si $j = i + ke$ et $k \le m_i$,} \\
\varnothing & \text{sinon.}
\end{cases}
\]
\end{cor}

\begin{proof}[D\'emonstration de la Proposition~\ref{prop:Gef1n-j-ind}]
Posons $g = fn$, et soit $\bgamma \in \cP(g,1,n)$ la multipartition telle que 
$j_{G(e,1,n)}^{G(g,1,n)} E_\balpha=E_\bgamma$. Par transitivit\'e de 
l'induction tronqu\'ee, on sait que $j_{G(ef,1,n)}^{G(g,1,n)} E_\bbeta 
\simeq E_\bgamma$ aussi.

Le Corollaire~\ref{cor:Gef1n-j-ind} d\'ecrit $\bgamma$ soit en termes de 
$\balpha$, soit en termes de $\bbeta$.

Puisque toute repr\'esentation irr\'eductible de $G(e,1,n)$ s'obtient par 
induction tronqu\'ee d'un produit tensoriel externe de repr\'esentations de la 
forme $\epsilon \otimes \gamma_e^k$, on peut d\'eduire du lemme pr\'ec\'edent une formule pour $\bgamma$ en fonction de $\balpha$. 
\end{proof}

La proposition suivante d\'ecoule imm\'ediatement de la
Proposition~\ref{prop:Gef1n-j-ind}.

\begin{prop}\label{prop:sym-00-dist}
Soit $\bbeta$ la multipartition d\'efinie par $E_\bbeta =
j_{\fS_n}^{G(e,1,n)} \epsilon$. Alors $\Lambda^{0,0}_\bb(\bbeta)$ est
distingu\'e.
\end{prop}

\begin{proof}[D\'emonstration de la Proposition~\ref{prop:paraboliques}]
Les sous-groupes paraboliques de $G(e,1,n)$ sont tous de la forme
\[
W' = G(e,1,n_0) \times \fS_{n_1} \times \cdots \times \fS_{n_k}
\qquad\text{o\`u}\qquad
n_0 + n_1 + \cdots + n_k = n.
\]
Soit $E$ une repr\'esentation sp\'eciale de $W'$: donc
\[
E = E_\balpha \boxtimes F_1 \boxtimes \cdots \boxtimes F_k,
\]
o\`u $\Lambda^{1,0}_\bb(\balpha)$ est distingu\'e, et $F_i$ est une repr\'esentation quelconque de $\fS_i$ (toutes ses repr\'esentations \'etant sp\'eciales). Rappelons que toute repr\'esentation du groupe sym\'etrique est l'induite tronqu\'ee de la repr\'esentation signe d'un produit de groupes sym\'etriques plus petits. Nous nous int\'eressons \`a $j_{W'}^{G(e,1,n)} E$, et donc on peut sans perte de g\'en\'eralit\'e supposer que $F_i = \epsilon$ pour tout $i$.

Soit $\bbeta_i$ la multipartition de la repr\'esentation $j_{\fS_{n_i}}^{G(e,1,n_i)} \epsilon$. Selon la Proposition~\ref{prop:sym-00-dist}, $\Lambda^{0,0}_\bb(\bbeta_i)$ est distingu\'e pour tout $i$. Il s'ensuit que la multipartition
\[
\balpha + \bbeta_1 + \cdots + \bbeta_k,
\]
qui correspond \`a $j_{W'}^{G(e,1,n)} E$, a la propri\'et\'e que son symbole de type $(1,0)$ est distingu\'e. Autrement dit, la repr\'esentation
$j_{W'}^{G(e,1,n)} E$ est sp\'eciale.

La preuve pour $G(e,e,n)$ est presque la m\^eme; on utilise le
Corollaire~\ref{cor:Geen-Ge1n} pour \'etudier les multipartitions des
$j_{\fS_{n_i}}^{G(e,e,n_i)} \epsilon$. (Il est \`a noter que cette
derni\`ere re\-pr\'esentation n'est pas forc\'ement bien d\'efinie: il peut
exister plusieurs exemplaires non conjugu\'es de $\fS_{n_i}$ dans
$G(e,e,n_i)$. N\'eanmoins les diverses repr\'esentations ainsi obtenues sont
toutes associ\'ees \`a la m\^eme multipartition).
\end{proof}

%**************************************************************************
\section{Repr\'esentations de Springer}
\label{sect:springer}
%**************************************************************************

Soit $W$ le groupe de Weyl d'un groupe alg\'ebrique r\'eductif.
Une repr\'esentation irr\'eductible de $W$ est dit \emph{de Springer} si
elle correspond, via la correspondance de Springer, au syst\`eme local
trivial sur une classe unipotente. En particulier, l'ensemble de
repr\'esentations de Springer est en bijection avec l'ensemble de classes
unipotentes.

Dans cette section, nous commen\c cons par rappeler une caract\'erisation
bien connue des repr\'esentations de Springer (voir le
Th\'eor\`eme~\ref{thm:springer}) en termes des sous-groupes
\emph{pseudoparaboliques}. Ensuite, on voudrait prendre l'\'enonc\'e
de ce th\'eor\`eme comme d\'efinition dans le cadre des groupes de
r\'eflexions complexes, mais il faudra l\'eg\`erement modifier la
d\'efinition de ``pseudoparabolique'' pour que les re\-pr\'esentations de
Springer soient bien d\'efinies. Apr\`es avoir trouv\'e la bonne
d\'efinition, nous donnons la classification des sous-groupes
pseudoparaboliques ainsi que celle des re\-pr\'esentations de Springer.

Il est bien connu que les diverses conditions figurant dans la d\'efinition
suivante sont \'equivalentes. 

\begin{defn}
Soit $W$ le groupe de Weyl d'un groupe alg\'ebrique r\'eductif $G$ \`a tore maximal $T$. Soit $G^*$ le groupe dual de $G$, et soit $T^*$ le tore maximal dual \`a $T$. Un sous-groupe $W' \subset W$ est \emph{\-pseudoparabolique} (\`a l'\'egard de $G$) s'il v\'erifie l'une des conditions suivantes \'equivalentes:
\begin{enumerate}
\item $W'$ est le centralisateur d'un point de $T^*$.
\item $W'$ est conjugu\'e \`a un sous-groupe engendr\'e par les
r\'eflexions correspondant \`a un sous-ensemble propre des n\oe uds du
diagramme de Dynkin \'etendu de $G^*$.
\end{enumerate}
Dans le cas o\`u $G$ et $G^*$ sont d\'efinis sur $\C$, on peut ajouter une
troisi\`eme version: soit $T^*_0$ une forme r\'eelle compacte de $T^*$, et soit
$\gt^*_0$ son alg\`ebre de Lie. Soit $L \subset \gt^*_0$ le r\'eseau
radiciel associ\'e \`a $G$ et $T$. On peut identifier $T^*_0$ avec
$\gt^*_0/L$, et donc $W'$ est pseudoparabolique si
\begin{enumerate}
\setcounter{enumi}{2}
\item $W'$ est le centralisateur d'un point de $\gt^*_0/L$.
\end{enumerate}
\end{defn}

La deuxi\`eme condition ci-dessus est importante car elle rend \'evident le
fait que tout sous-groupe pseudoparabolique est engendr\'e par des
r\'eflexions; pourtant, dans le cadre des groupes de r\'eflexions
complexes, o\`u on ne dispose pas d'une th\'eorie bien d\'evelopp\'ee de
diagrammes de Dynkin, c'est la troisi\`eme condition qui pourra se
g\'en\'eraliser. La caract\'erisation suivante des repr\'esentations de 
Springer est bien connue (voir par exemple \cite[\S 12.6]{Ca}).

\begin{thm}\label{thm:springer}
Soit $W$ le groupe de Weyl d'un groupe alg\'ebrique r\'eductif $G$. Une
re\-pr\'esentation
de $W$ est de Springer si et seulement si elle est l'induite tronqu\'ee
d'une repr\'esentation sp\'eciale d'un sous-groupe pseudoparabolique. 
\end{thm}

D\'esormais, nous travaillons dans le contexte suivant: $K$ d\'esigne un corps
de nombres ab\'elien, $V$ est un espace vectoriel de dimension finie
sur $K$, $W \subset GL(V)$ est un groupe de r\'eflexions et $V$ est muni
d'une forme hermitienne qui est $W$-invariante et non d\'eg\'en\'er\'ee. Enfin, soit $\Z_K$ l'anneau des entiers alg\'ebriques dans $K$.

\begin{defn}[Nebe]
Une \emph{racine} pour
$W$ est un vecteur propre pour une r\'eflexion dans $W$ \`a valeur propre
non triviale. 

Un \emph{r\'eseau radiciel primitif} pour $W$ est un $\Z_K$-sous-module de
$V$ qui est $W$-invariant et engendr\'e en tant que $\Z_K[W]$-module par une 
seule racine.
\end{defn}

Nebe a classifi\'e dans \cite{N} tous les r\'eseaux radiciels primitifs des 
groupes de r\'eflexions complexes. Plus tard, nous rappellerons ses r\'esultats 
pour les groupes imprimitifs sp\'etsiaux, mais d'abord, esquissons ce que nous
voudrions faire:
\begin{enumerate}
\item D\'efinir les sous-groupes pseudoparaboliques de $W$, \`a l'\'egard
d'un r\'eseau radiciel $L$, comme \'etant les stabilisateurs des points de
$V/L$.
\item D\'efinir les repr\'esentations de Springer de $W$ \`a l'\'egard de
$L$ comme \'etant les induites tronqu\'ees des repr\'esentations
sp\'eciales des sous-groupes pseudoparaboliques.
\item Esp\'erer que tous les sous-groupes pseudoparaboliques aient la forme
des sous-groupes figurant dans les Th\'eor\`emes~\ref{thm:Ge1n-springer}
et~\ref{thm:Geen-springer}, et puis d\'eduire que les symboles distingu\'es
param\`etrent les repr\'esentations de Springer. 
\end{enumerate}
Malheureusement, ce projet \'echoue \`a la premi\`ere \'etape: les
stabilisateurs des points de $V/L$ ne sont m\^eme pas toujours des groupes
de r\'eflexions, et si l'on se restreint aux sous-groupes du stabilisateur
engendr\'es par des r\'eflexions, on peut obtenir des groupes qui sont non
sp\'etsiaux ou non \emph{pleins} (voir la D\'efinition~\ref{defn:plein}).

Nous allons bien associer \`a chaque point de $V/L$ un sous-groupe de $W$
(ou plut\^ot, une classe de conjuaison de sous-groupes de $W$) qui sera dit
``pseudoparabolique'', et puis nous effectuerons le reste de l'esquisse comme
d\'ecrite ci-dessus. Pourtant, le besoin d'\'eviter tous ces probl\`emes complique beaucoup la construction.

Si $L$ est un r\'eseau radiciel primitif pour $W$, soit $\tWL$ le groupe
engendr\'e par $W$ et toute autre r\'eflexion qui pr\'eserve $L$ et
pour laquelle il y a une racine dans $L$. Alors $\tWL$ est un groupe
de r\'eflexions contenant $W$ (et \'eventuellement plus grand) pour lequel
$L$ est un r\'eseau radiciel primitif.

Soit $x \in V/L$. Soit $(\tWL)^x$ le stabilisateur dans $\tWL$ de $x$. Ce
groupe-ci n'est pas forc\'ement engendr\'e par des r\'eflexions, et donc on
d\'efinit $(\tWL)^x_{\refl}$ comme \'etant le groupe engendr\'e par les
r\'eflexions dans $(\tWL)^x$. Maintenant, il n'est pas forc\'ement le cas
que $(\tWL)^x_{\refl}$ ait la propri\'et\'e suivante, laquelle sera
essentielle pour les sous-groupes pseudoparaboliques.

\begin{defn}\label{defn:plein}
Un sous-groupe de r\'eflexions $W' \subset W$ est \emph{plein} si toute
r\'eflexion $s \in W$ v\'erifie la propri\'et\'e suivante: s'il y a un
entier $a$ tel que $s^a \in W'$ et $s^a \ne 1$, alors $s \in W'$. 
\end{defn}

Par un abus de notation, nous d\'efinissons $(\tWL)^x_{\plein}$ comme
\'etant un sous-groupe parabolique maximal de $(\tWL)^x_{\refl}$ qui est
plein en tant que sous-groupe de $\tWL$. (C'est un abus car
$(\tWL)^x_{\plein}$ n'est unique qu'\`a conjugaison pr\`es). Enfin,
$(\tWL)^x_{\plein}$ n'est pas forc\'ement sp\'etsial. On note
$(\tWL)^x_{\spets}$ le plus grand sous-groupe sp\'etsial plein de
$(\tWL)^x_{\plein}$.

\begin{defn}\label{defn:pseudopar}
Soit $W$ un groupe de r\'eflexions sur $V$, et soit $L$ un r\'eseau
radiciel primitif pour $W$. Un sous-groupe $W' \subset W$ est un
sous-groupe \emph{pseudoparabolique associ\'e \`a $x \in V/L$} si $W'
= W \cap (\tWL)^x_{\spets}$ pour un certain groupe $(\tWL)^x_{\spets}$.

Une repr\'esentation irr\'eductible de $W$ est dite \emph{de Springer} si
elle est l'induite tronqu\'ee d'une repr\'esentation sp\'eciale d'un
sous-groupe pseudoparabolique.
\end{defn}

Nous allons d\'eterminer tous les sous-groupes pseudoparaboliques ainsi
que les repr\'esentations de Springer pour tous les groupes de
r\'eflexions complexes imprimitifs. Il faut traiter les groupes di\'edraux
s\'epar\'ement car leurs corps de d\'efinition sont de la forme $\Q(\zeta
+ \zeta^{-1})$ (o\`u $\zeta$ est une racine de l'unit\'e) plut\^ot que
$\Q(\zeta)$, et par cons\'equent, la d\'etermination des groupes $\tilde
W_L$ ainsi que celle des sous-groupes pleins ou sp\'etsiaux est
diff\'erente.

%--------------------------------------------------------------------------
\subsection{Les groupes imprimitifs non di\'edraux}
%--------------------------------------------------------------------------

Soient $e$ un entier positif et $\zeta$ une racine $e$-\`eme
de l'unit\'e primitive, et posons $K = \Q(\zeta)$, $V = \bigoplus_{i=1}^n Ke_i$ ({\it cf.} la Section~\ref{subsection: repr}), $L_1 =
\bigoplus_{i=1}^n \Z_Ke_i$, et
\[
L_2 = \{\textstyle\sum v_i e_i \in L_1 \mid \sum v_i \in
(1-\zeta)\Z_K \}.
\]
(Il est \`a noter que si $e$ n'est pas une puissance d'un nombre premier,
alors $1-\zeta$ est inversible dans $\Z_K$, et donc $L_1 = L_2$). On
pose aussi $W_n = G(e,1,n)$ et $W'_n = G(e,e,n)$.

\begin{thm}[Nebe]\label{thm:nebe}
\begin{enumerate}
\item Si $e$ n'est pas une puissance d'un nombre premier, alors
$L_1$ est l'unique r\'eseau radiciel primitif de $W_n$ \`a
isomorphisme pr\`es.
\item Si $e$ est une puissance d'un nombre premier, alors $W_n$ admet
deux classes d'isomorphie de r\'eseaux radiciels primitifs, dont $L_1$ et
$L_2$ sont des re\-pr\'esentants.
\item $L_2$ est l'unique r\'eseau radiciel primitif de $W'_n$ \`a
isomorphisme pr\`es si $n \ge 3$.
\end{enumerate}
\end{thm}

\`A la suite de ce th\'eor\`eme et la classification des groupes de
r\'eflexions complexes, l'\'enonc\'e suivant est \'evident.

\begin{lem}
Pour tout groupe de r\'eflexions complexes imprimitif irr\'eductible $W$
d\'efini sur $K$, et tout r\'eseau radiciel primitif $L$ pour $W$, on a
que $\tilde W_L = G(e,1,n)$.
\end{lem}

\begin{lem}\label{lem:bloc-stab}
Soit $v = \sum v_i e_i \in V$ un point tel que $v_i \equiv v_j \pmod
{\Z_K}$ pour tout $i$, $j$. Soit $x$ l'image de $v$ dans $V/L$, o\`u $L$ est un r\'eseau radiciel primitif. Si $v_i \in \Z_K$ pour tout $i$, alors
\[
(\tilde W_{L})^x_{\refl} =
(\tilde W_{L})^x_{\plein} =
(\tilde W_{L})^x_{\spets} =
G(e,1,n).
\]
Sinon, soit $t$ le plus petit entier
strictement positif tel
que $(1-\zeta^t)v_1 \in \Z_K$ (et donc $(1-\zeta^t)v_i \in \Z_K$ pour tout
$i$). On a que
\begin{align*}
(\tilde W_{L_1})^x_{\refl} &= G(e/t, 1, n),
%\begin{cases}
%G(e/t, 1, n) & \text{si $t < e$,} \\
%G(1,1,n) & \text{si $t = e$,}
%\end{cases}
\\
(\tilde W_{L_1})^x_{\plein} = (\tilde W_{L_1})^x_{\spets} &=
\begin{cases}
G(e,1,n) & \text{si $t = 1$,} \\
G(1,1,n) & \text{si $t > 1$.}
\end{cases}
\end{align*}

Si $e$ est une puissance d'un nombre premier $p$, alors on a
aussi
\begin{align*}
(\tilde W_{L_2})^x_{\refl} &=
\begin{cases}
G(e/t, p, n) & \text{si $t < e$,} \\
G(1,1,n) & \text{si $t = e$,}
\end{cases}
\\
(\tilde W_{L_2})^x_{\plein} &=
\begin{cases}
G(e,p,n) \\%& \text{si $t = 1$,} \\
G(1,1,n) %& \text{si $t > 1$,}
\end{cases}
\quad\text{et}\quad
(\tilde W_{L_2})^x_{\spets} =
\begin{cases}
G(e,e,n) & \text{si $t = 1$,} \\
G(1,1,n) & \text{si $t > 1$.}
\end{cases}
\end{align*}
\end{lem}
\begin{proof}
Le cas o\`u les $v_i$ sont dans $\Z_K$ est \'evident. Nous supposons d\'esor\-mais que les $v_i \notin \Z_K$. Remarquons qu'une fois qu'on connait les $(\tilde W_{L})^x_{\refl}$, les $(\tilde W_{L})^x_{\plein}$ et les $(\tilde W_{L})^x_{\spets}$ s'en d\'eduisent tr\`es facilement. Il suffit de d\'eterminer les $(\tilde W_{L})^x_{\refl}$.  Remarquons aussi que $t$ divise $e$.  En particulier, si $e$ est une puissance d'un nombre premier, alors $t$ l'est aussi.

Soit $s \in \fS_n \subset G(e,1,n)$ la transposition qui \'echange $e_1$ et $e_2$. Puisqu'on a suppos\'e que $v_i \equiv v_j \pmod{\Z_K}$, on voit que
\[
v - s\cdot v = (v_1 - v_2)e_1 + (v_2 - v_1)e_2 \in L_2 \subset L_1.
\]
On peut faire un calcul analogue pour toute transposition, et donc on sait
que $\fS_n$ stabilise $x$. Ensuite, soit $r \in G(e,1,n)$ la r\'eflexion
qui envoie $e_1$ sur $\zeta e_1$ et fixe les autres $e_i$.  Alors
\begin{multline*}
v - (r^{-t} s r^t)v = (v_1 - \zeta^{-t}v_2) e_1 + (v_2 - \zeta^t v_1) e_2 \\
= ((v_1 - v_2) -\zeta^{-t}(1-\zeta^t)v_2)e_1 + ((v_2-v_1) + (1-\zeta^t)v_1)e_2 \in L_1.
\end{multline*}
Soit $n = v_2 - v_1 \in \Z_K$. Alors la somme des coefficients de cette expression \'egale
\[
-\zeta^{-t}(1-\zeta^t)v_2 + (1-\zeta^t)(v_2 - n)
= -\zeta^{-t}(1-\zeta^t)(1 - \zeta^t)v_2 + (1-\zeta^t)n  \in (1-\zeta^t)\Z_K.
\]
On voit que $v - (r^{-t} s r^t)v \in L_2$. En combinaison avec $\fS_n$, l'\'el\'ement $r^{-t}sr^t$ engendre le groupe $G(e/t,e/t,n)$. Donc $G(e/t, e/t, n) \subset (\tilde W_{L})^x_{\refl}$ et pour $L = L_1$, et pour $L = L_2$.

Supposons maintenant que $L = L_1$, et soit $k$ un facteur de $e$.  Il est clair que
\[
v - r^k v= (1- \zeta^k)v_1e_1 \quad
\begin{cases}
\in L_1 & \text{si $k = t$,} \\
\notin L_1 & \text{si $1 \le k < t$.}
\end{cases}
\]
Donc $r^t \in (\tilde W_{L_1})^x_{\refl}$, mais $r^k \notin (\tilde W_{L_1})^x_{\refl}$ si $1 \le k < t$.  On voit  que $G(e/t,1,n) \subset (\tilde W_{L_1})^x_{\refl}$; de surcro\^it, puisqu'on a d\'ej\`a consid\'er\'e toutes les r\'eflexions dans $\tilde  W_{L_1}$ (\`a conjugaison pr\`es), on conclut que $(\tilde W_{L_1})^x_{\refl} = G(e/t,1,n)$.

Supposons maintenant que $e$ est une puissance d'un nombre premier. Si $t = e$, alors aucune puissance non triviale de $r$ ne fixe $x$, et $(\tilde W_{L_2})^x_{\refl} = G(1,1,n)$. Supposons d\'esormais que $t < e$. M\^eme si $k = t$, il n'est pas vrai que $v - r^k v \in L_2$, car $(1-\zeta^t)v_1 \in \Z_K$ mais $(1-\zeta^t)v_1 \notin (1-\zeta)\Z_K$. Rappelons que $(1-\zeta^{tp})\Z_K = (1-\zeta^t)^p\Z_K$. Par cons\'equent,
\[
v - r^{tp} v = (1-\zeta^{tp})v_1 e_1 \in (1-\zeta^t)^p \Z_K v_1
\subset (1-\zeta^t)^{p-1} \Z_K e_1 \subset (1-\zeta)^{p-2} L_2.
\]
$(\tilde W_{L_2})^x_{\refl}$ est engendr\'e par $G(e/t,e/t,n)$ et $r^{tp}$, et donc on trouve que $(\tilde W_{L_2})^x_{\refl} = G(e/t, p, n)$.
\end{proof}

Enfin, nous pouvons d\'emontrer le th\'eor\`eme principal.

\begin{thm}\label{thm:princ}
Soit $W$ un groupe de r\'eflexions complexes imprimitif irr\'eductible
sp\'etsial non di\'edral, et soit $L$ un r\'eseau radiciel primitif pour
$W$. L'ensemble
des sous-groupes pseudoparaboliques (\`a conjugaison pr\`es) de $W$ \`a
l'\'egard de $L$ est indiqu\'e ci-dessous.

Une repr\'esentation de $W$ est
de Springer \`a l'\'egard de $L$ si et seulement si son symbole de poids
sp\'etsial et d'un certain type (qui d\'epend de $W$ et de $L$) est
distingu\'e. La table ci-dessous pr\'ecise le type convenable pour chaque
groupe $W$ et chaque r\'eseau $L$.
\[
\begin{array}{c|l|c}
\text{groupe;} &\multicolumn{1}{c|}{\text{sous-groupes}} & \text{type
des} \\
\text{r\'eseau} & \multicolumn{1}{c|}{\text{pseudoparaboliques}} &
\text{symboles} \\
\hline\hline
\multicolumn{3}{l}{\text{$G(e,1,n)^{\strut}$,\quad o\`u $e = p^a$,
$p$ premier}} \\
\hline
L_1 & \prod_{i=0}^{p-1} G(e,1,n_i) \times \prod_{i=1}^k \fS_{m_i} & (p,0)
\\
L_2& G(e,1,n_0) \times \prod_{i=1}^{p-1} G(e,e,n_i) \times \prod_{i=1}^k
\fS_{m_i} & (p,p-1) \\
\hline\hline
\multicolumn{3}{l}{\text{$G(e,1,n)^{\strut}$,\quad o\`u $e$ n'est pas une
puissance d'un nombre premier}} \\
\hline
L_1=L_2 & G(e,1,n_0) \times \prod_{i=1}^k \fS_{m_i} & (1,0) \\
\hline\hline
\multicolumn{3}{l}{\text{$G(e,e,n)^{\strut}$,\quad o\`u $e = p^a$, $p$
premier, et $n \ge 3$}} \\
\hline
L_2 & \prod_{i=0}^{p-1} G(e,e,n_i) \times \prod_{i=1}^k \fS_{m_i} & (p,0)
\\
\hline\hline
\multicolumn{3}{l}{\text{$G(e,e,n)^{\strut}$,\quad o\`u $e$ n'est pas une
puissance d'un nombre premier, et $n \ge 3$}} \\
\hline
L_1=L_2 & G(e,e,n_0) \times \prod_{i=1}^k \fS_{m_i} & (1,0) \\
\hline
\end{array}
\]
Ici, on a que $k \ge 0$, et
les $n_i$ et les $m_i$ sont des entiers tels $n_i \ge 0$, $m_i \ge 1$, et $\sum n_i + \sum m_i = n$.
\end{thm}

\begin{rem}
Le param\'etrage des repr\'esentations de Springer fourni par ce th\'eor\`eme pour $G(2,1,n)$ et $L_1$ (resp.~$G(2,1,n)$ et $L_2$, $G(2,2,n)$) co\"incide avec celui provenant des classes unipotentes et de la correspondance de Springer pour un groupe alg\'ebrique de type $B_n$ (resp.~$C_n$, $D_n$).
\end{rem}

\begin{proof}
Remarquons d'abord il suffit de calculer les sous-groupes pseudoparaboliques dans chaque cas; le type convenable des symboles se d\'eduit des r\'esultats des Sections~\ref{sect:Ge1n} et~\ref{sect:paraboliques}. Soit $x \in V/L$, et soit $v \in V$ un point dans l'image r\'eciproque de $x$. On peut \'evidemment remplacer $x$ et $v$ par d'autres points dans leurs $G(e,1,n)$-orbites respectives sans changer la classe d'isomorphie du sous-groupe pseudo\-parabolique associ\'e. On peut donc imposer l'hypoth\`ese suivante sur $v = (v_1, \ldots, v_n)$ sans perte de g\'en\'eralit\'e: les coordonn\'ees $v_1, \ldots, v_n$ se r\'epartissent en ``blocs''
\[
v_1, \ldots, v_{a_1};\quad v_{a_1+1}, \ldots, v_{a_2};\quad \ldots;\quad
v_{a_{l-1}+1}, \ldots, v_{a_l}
\qquad\text{(o\`u $a_l = n$)}
\]
tel que
\[
\begin{aligned}
v_i & \equiv v_j \pmod{\Z_K} &&\text{si $i$ et $j$ appartiennent au m\^eme
bloc,} \\
v_i & \not\equiv \zeta^k v_j \pmod{\Z_K} &&\text{pour tout $k$ sinon.}
\end{aligned}
\]

Si $e$ est une puissance d'un nombre premier, soit $p$ ledit nombre
premier; sinon, posons $p = 1$. Rappelons que $1-\zeta$ engendre un id\'eal maximal de $\Z_K$ au-dessus de $(p) \subset \Z$ si $e$ est une puissance d'un nombre premier, et est inversible sinon. Dans tous les deux cas, on a que $\Z_K/(1-\zeta)\Z_K \simeq (1-\zeta)^{-1}\Z_K/\Z_K \simeq \Z/p\Z$. \'Evidemment, les entiers $0,1,\ldots, p-1$ d\'ecrivent un ensemble de repr\'esentants des \'elements de $\Z_K/(1-\zeta)\Z_K$ ou de $\Z/p\Z$.  De m\^eme, les \'elements $k/(1-\zeta)$, o\`u $k \in \{0,1 \ldots, p-1\}$, d\'ecrivent un ensemble de repr\'esentants de $(1-\zeta)^{-1}\Z_K/\Z_K$. En particulier, il y a au plus $p$ blocs de $v$ dont les membres appartiennent \`a $(1-\zeta)^{-1}\Z_K$. Supposons, sans perte de g\'en\'eralit\'e, que ces blocs-l\`a soient les $p$ premiers blocs, et que
\[
v_i \equiv \frac{k}{1-\zeta} \pmod{\Z_K}
\qquad\text{si $k \in \{0,1,\ldots,p-1\}$ et $a_{k}+1 \le i \le a_{k+1}$.}
\]
(On permet que $a_{k-1} = a_k$, \ie que certains blocs soient vides). Pour s'assurer que les blocs ainsi d\'efinis sont bien distincts, il faut d\'emontrer que $k/(1-\zeta) \not\equiv \zeta^l k'/(1-\zeta) \pmod{\Z_K}$ pour tout $l$ si $k, k' \in \{0, \ldots, p-1\}$ et $k \ne k'$. Ceci est tr\`es facile:
\begin{align*}
\frac{k}{1-\zeta} - \frac{\zeta^l k'}{1-\zeta} &=
\frac{k - k'}{1-\zeta} + \frac{k' - \zeta^l k'}{1-\zeta} 
= \frac{k-k'}{1-\zeta} + (1 + \zeta + \cdots + \zeta^{l-1})k' \\
&\equiv \frac{k-k'}{1-\zeta} \pmod{\Z_K} \not\equiv 0 \pmod{\Z_K}.
\end{align*}

On peut maintenant invoquer le Lemme~\ref{lem:bloc-stab} pour chaque bloc.
Dans le premier bloc, on a $v_1, \ldots, v_{a_1} \in \Z_K$, et le calcul
du groupe $(\tilde W_{m,L})^x_{\spets}$ l\`a-bas donne le premier facteur
de chaque sous-groupe pseudoparabolique dans la table ci-dessus.

Ensuite, consid\'erons un bloc $v_{a_k+1}, \ldots, v_{a_{k+1}}$ o\`u $1
\le k < p$. Dans le cadre du Lemme~\ref{lem:bloc-stab}, on a $t = 1$, et on
obtient ainsi $p-1$ facteurs de type $G(e,1,m)$ ou $G(e,e,m)$, selon le
r\'eseau, dans chaque sous-groupe pseudoparabolique.  Enfin, dans tous les
bloc apr\`es le $p$-\`eme, on invoque ce lemme-l\`a avec $t > 1$ , et on
trouve que $(\tilde W_{m,L})^x_{\spets}$ est toujours un groupe
sym\'etrique.
\end{proof}

%--------------------------------------------------------------------------
\subsection{Les groupes di\'edraux}
\label{sect:diedraux}
%--------------------------------------------------------------------------

On garde les notations de la section pr\'ec\'edente: $e$ est un entier
positif, $\zeta$ est une $e$-\`eme racine de l'unit\'e primitive, $K =
\Q(\zeta)$ et $\Z_K = \Z[\zeta]$. Soit $V = Ke_1 \oplus K e_2$.  Nous
consid\'erons le groupe $W = G(e,e,2)$. Le corps de d\'efinition de $W$
est $K_0 = \Q(\zeta + \zeta^{-1})$, mais pour certains calculs
ult\'erieurs, il sera commode d'avoir d\'efini $W$ sur $K$.

Rappelons que tout sous-groupe de r\'eflexions de $W$ est isomorphe \`a
$G(d,d,2)$, o\`u $d \mid e$. Posons 
\[
s_i = \begin{bmatrix} 0 & \zeta^i \\ \zeta^{-i} & 0 \end{bmatrix},
\]
et pour tout diviseur $d$ de $e$ ($e$ compris), identifions $G(d,d,2)$ avec
le sous-groupe de $GL(V)$ engendr\'e par $s_0$ et $s_{e/d}$. Dans le cas
o\`u $e/d$ est pair, on note $G'(d,d,2)$ le sous-groupe engendr\'e par
$s_1$ et $s_{e/d+1}$. Ce dernier est isomorphe mais non conjugu\'e \`a
$G(d,d,2)$. Tout sous-groupe de r\'eflexions de $W$ est conjugu\'e ou bien
\`a l'un des $G(d,d,2)$ ou bien \`a l'un des $G'(d,d,2)$.

Rappelons la classification des repr\'esentations irr\'eductibles des
groupes di\'edraux. Celles du groupe $G(d,d,2)$
seront not\'ees: 
\[
\chi^{(d)}_0, \chi^{(d)}_1, \ldots, \chi^{(d)}_{\lfloor(d-1)/2\rfloor};
\qquad \chi^{(d)}_d;
\qquad\text{et, si $d$ est pair,} \qquad
\chi^{(d)}_{d/2}, \chi^{(d)\prime}_{d/2}.
\]
(Les indices en bas indiquent les valeurs de la fonction $b$). La
repr\'esentation triviale ($\chi^{(d)}_0$), la repr\'esentation r\'eflexion
($\chi^{(d)}_1$), et la repr\'esentation signe ($\chi^{(d)}_d$) sont les
seules repr\'esentations sp\'eciales. 

L'induite tronqu\'ee de la repr\'esentation triviale
(resp.~r\'e\-flex\-ion) \`a partir de n'importe quel $G(d,d,2)$ ou
$G'(d,d,2)$ \`a $G(e,e,2)$ est encore la repr\'esentation triviale
(resp.~r\'eflexion). Par contre, on a: 
\[
j_{G(d,d,2)}^{G(e,e,2)} \chi^{(d)}_d = \chi^{(e)}_d
\qquad\text{et}\qquad
j_{G'(d,d,2)}^{G(e,e,2)} \chi^{(d)}_d =
\begin{cases}
\chi^{(e)}_d & \text{si $d \ne e/2$,} \\
\chi^{(e)\prime}_{e/2} & \text{si $d = e/2$.}
\end{cases}
\]

Enfin, nous avons besoin d'un r\'eseau radiciel pour $W$. Posons
\[
V_0 = \{ v_1 e_1 + v_2 e_2 \in V \mid v_2 = - \bar v_1\}.
\]
Il est facile de v\'erifier que $V_0$ est un $K_0$-sous-espace de $V$ qui est stable sous $G(e,e,2)$, et que $V \simeq V_0 \otimes_{K_0} K$. Ensuite, posons
\[
L_0 = \{ v_1 e_1 - \bar v_1 e_2 \in V_0 \mid v_1 \in \Z_K \}.
\]
$L_0$ est un $\Z_{K_0}$-r\'eseau dans $V_0$, stable sous $G(e,e,2)$ et
engendr\'e par la $G(e,e,2)$-orbite du vecteur $e_1-e_2 \in V_0$. Ce
dernier \'etant une racine pour $s_0$, on voit que $L_0$ est bien un
r\'eseau radiciel primitif. Rappelons maintenant le r\'esultat de Nebe sur
les r\'eseaux radiciels des groupes di\'edraux.

\begin{thm}[Nebe]\label{thm:nebe-diedr}
$L_0$ est un repr\'esentant de l'unique genre de r\'eseaux radiciels
primitifs sauf si $e$ est pair et $e/2$ est une puissance d'un nombre
premier. Dans ce dernier cas, il y deux genres de r\'eseaux radiciels
primitifs, dont l'un est repr\'esent\'e par $L_0$, et l'autre par le
r\'eseau engendr\'e par la $W$-orbite d'une racine pour $s_1$.
\end{thm}

Pourtant, dans le cas o\`u $e$ est pair et $e/2$ est une
puissance d'un nombre premier, l'automorphisme externe qui \'echange les
deux classes de conjugaison de r\'eflexions \'echange aussi les deux genres
de r\'eseaux radiciels ({\it cf.} les r\'eseaux radiciels des groupes de
Weyl de type $B_2 = G(4,4,2)$ ou $G_2 = G(6,6,2)$). Pour le calcul des
sous-groupes pseudoparaboliques, il suffit de consid\'erer
seulement le r\'eseau $L_0$.

L'analogue du Lemme~\ref{lem:bloc-stab} est \'evident:

\begin{lem}
Pour le groupe groupe di\'edral $W = G(e,e,2)$, on a que $\tilde W_{L_0} =
W$.
\end{lem}

\begin{thm}\label{thm:princ-diedr}
Le sous-groupe $G(d,d,2)$ de $G(e,e,2)$ est pseudoparabolique si et seulement si $d$ v\'erifie l'une des conditions suivantes:
\begin{itemize}
\item $d = 1$,
\item $d = e$,
\item $d$ divise $e$ et est une puissance d'un nombre premier.
\end{itemize}
Si $e$ et pair, $G'(d,d,2)$ est pseudoparabolique si et seulement si $e/d$ est un entier pair et l'une des conditions suivantes est satisfaite:
\begin{itemize}
\item $d = 1$,
\item $d < e/2$ et $d$ est une puissance d'un nombre premier.
\end{itemize}
Si $e > 2$, l'ensemble de repr\'esentations de Springer de $G(e,e,2)$ est:
\[
\{\chi^{(e)}_0, \chi^{(e)}_1, \chi^{(e)}_e\} \cup
\{\chi^{(e)}_d \mid \text{$d$ divise $e$ et est une puissance d'un nombre premier} \}.
\]
En particulier, si $e$ est pair et sup\'erieur \`a $2$, $\chi^{(e)\prime}_{e/2}$ n'est pas de Springer.

Toute repr\'esentation irr\'eductible de $G(2,2,2)$ est de Springer.
\end{thm}

\begin{exm}
Les sous-groupes pseudoparaboliques de $G(6,6,2) = G_2$ sont de type
$G(1,1,2) = A_1$, $G(2,2,2) = A_1 \times A_1$, et $G(3,3,2) = A_2$. Il y a
aussi un deuxi\`eme exemplaire de $A_1 = G'(1,1,2)$, qui est lui aussi
pseudoparabolique, ainsi que de $A_2 = G'(3,3,2)$, qui ne l'est pas. Les
repr\'esentations de Springer de $G_2$ sont $\chi_0$, $\chi_1$, $\chi_2$,
$\chi_3$, et $\chi_6$ ({\it cf.}~\cite{Ca}).
\end{exm}

\begin{proof}
Le cas du groupe $G(2,2,2)$ est tr\`es facile: ce groupe-l\`a est isomorphe \`a $\fS_2 \times \fS_2$, dont toute repr\'esentation est sp\'eciale et donc de Springer. On suppose d\'esormais que $e > 2$.

La liste des repr\'esentations de Springer se d\'eduit tr\`es facilement de
la liste des sous-groupes pseudoparaboliques et des rappels ci-dessus sur
l'induction tronqu\'ee. De plus, il est clair que tout groupe
$(\tilde W_{L_0})^x_{\refl} = W^x_{\refl}$, \'etant lui aussi un groupe
di\'edral, est d\'ej\`a plein et sp\'etsial. Il suffit donc de trouver les
$W^x_{\refl}$.

Soit $v = v_1 e_1 - \bar v_1 e_2 \in V_0$, et soit $x$ son image dans
$V_0/L_0$. Pour qu'une r\'eflexion $s_i$ appartienne \`a
$W^x_{\refl}$, il faut que $v - s_i v$ soit dans $L_0$. On a:
\begin{align*}
s_i v &=  -\zeta^i \bar v_1 e_1 + \zeta^{-i} v_1 e_2, \\
v - s_i v  &= (v_1 + \zeta^i \bar v_1) e_1 - (\bar v_1 + \zeta^{-i} v_1),
\end{align*}
et donc on voit que $s_i \in W^x_{\refl}$ si et seulement si
\begin{equation}\label{eqn:x-stab}
v_1 + \zeta^i \bar v_1 \in \Z_K.
\end{equation}

D\'emontrons d'abord que pour que $G(d,d,2)$ ou $G'(d,d,2)$ soit
pseudoparaboli\-que, il faut que $d$ soit $1$ ou une puissance d'un nombre
premier. Si $s_0$ et $s_{e/d}$ (resp.~$s_1$ et $s_{e/d+1}$) appartiennent
\`a $W^x_{\refl}$, alors 
\[
\begin{aligned}
v_1 + \bar v_1 &\in \Z_K \\
v_1 + \zeta^{e/d}\bar v_1 &\in \Z_K
\end{aligned}
\qquad\text{resp.}\qquad
\begin{aligned}
v_1 + \zeta \bar v_1 &\in \Z_K \\
v_1 + \zeta^{e/d+1}\bar v_1 &\in \Z_K
\end{aligned}
\]
et donc
\begin{equation}\label{eqn:stable}
(1-\zeta^{e/d})\bar v_1 \in \Z_K
\qquad\text{resp.}\qquad
\zeta(1-\zeta^{e/d})\bar v_1 \in \Z_K.
\end{equation}
Si $d$ n'est pas \'egal \`a $1$ ou une puissance d'un nombre premier, alors $1-\zeta^{e/d}$
est inversible dans $\Z_K$ ($\zeta^{e/d}$ \'etant une $d$-\`eme racine de
l'unit\'e primitive), et donc on voit que $\bar v_1 \in \Z_K$. Il s'ensuit
que $v \in L_0$ et que $W^x_{\refl} = W$. 

Avant de faire le prochain pas, rappelons que $K = K_0[\zeta]$ est une extension de $K_0$ de degr\'e $2$, et que $\Z_K = \Z_{K_0}[\zeta]$. Tout \'el\'ement de $K$ s'\'ecrit $\alpha + \beta\zeta$, o\`u $\alpha, \beta \in K_0$, de mani\`ere unique, et un tel \'el\'ement appartient \`a $Z_K$ si et seulement si $\alpha, \beta \in \Z_{K_0}$.

Ensuite, d\'emontrons que $G'(e/2,e/2,2)$ n'est pas pseudoparabolique. Si
$x$ est stable sous $G'(e/2,e/2,2)$, alors, selon~\eqref{eqn:stable}, on a
que $(1 - \zeta^2)\bar v_1 \in \Z_K$. L'expression pour $(1-\zeta^2)\bar v_1$ telle que d\'ecrite au paragraphe pr\'ec\'edent est
\[
(v_1 + \bar v_1) + (-\zeta^{-1}v_1 - \zeta \bar v_1)\zeta = (1-\zeta^2)\bar v_1.
\]
En particulier, on voit que $v_1 + \bar v_1 \in \Z_{K_0} \subset \Z_K$. \`A
la suite de~\eqref{eqn:x-stab}, on voit que $s_0 \in W^x_{\refl}$.
Puisqu'on avait d\'ej\`a suppos\'e que $s_1 \in W^x_{\refl}$, il s'ensuit
que $W^x_{\refl} = W$. 

Enfin, si l'on n'est pas dans les deux cas pr\'ec\'edents, on peut
construire explicitement un vecteur $v$ tel que $W^x_{\refl}$ \'egale
$G(d,d,2)$ (resp.~$G'(d,d,2)$). Posons 
\[
v_1 = \frac{1}{1-\zeta^{-e/d}}
\qquad\text{resp.}\qquad
v_1 = \frac{1+ \zeta}{1-\zeta^{-e/d}}
\]
et $v = v_1e_1 -\bar v_1 e_2$. Pour montrer que $G(d,d,2) \subset W^x_{\refl}$ (resp.~$G'(d,d,2) \subset W^x_{\refl}$), il suffit de montrer que l'\'egalit\'e~\eqref{eqn:stable} est v\'erifi\'ee pour $i = 0, e/d$ (resp.~$i = 1, e/d+1$).
En effet, apr\`es des calculs tr\`es faciles, on trouve que
\[
\begin{aligned}
v_1 + \bar v_1 &= 1, \\
v_1 + \zeta^{e/d} \bar v_1 &= 0
\end{aligned}
\qquad\text{resp.}\qquad
\begin{aligned}
v_1 + \zeta \bar v_1 &= 1 + \zeta, \\
v_1 + \zeta^{e/d+1} \bar v_1 &= 0.
\end{aligned}
\]
Par exemple, la premi\`ere \'egalit\'e ci-dessus se montre comme suit:
\[
v_1 + \bar v_1 = \frac{1}{1-\zeta^{-e/d}} + \frac{1}{1-\zeta^{e/d}} 
= \frac{(1-\zeta^{e/d}) + (1-\zeta^{-e/d})}{(1-\zeta^{-e/d})(1-\zeta^{e/d})}
= \frac{2 - \zeta^{e/d} - \zeta^{-e/d}}{2 - \zeta^{e/d} - \zeta^{-e/d}}
= 1.
\]
Il reste de s'assurer que $W^x_{\refl}$ ne soit pas plus grand qu'on ait voulu. Si l'on avait $W^x_{\refl} = G(f,f,2)$ (resp.~$G'(f,f,2)$) avec $f > d$, alors,
selon~\eqref{eqn:stable}, on aurait que $(1-\zeta^{e/f})\bar v_1 \in \Z_K$,
o\`u $e/f < e/d$. 

Par contre, nous verrons maintenant que le plus petit entier strictement
positif $t$ tel que $(1-\zeta^t)\bar v_1 \in \Z_K$ est $t = e/d$. C'est
clair dans le cas de $G(d,d,2)$, o\`u on a $\bar v_1 = 1/(1-\zeta^{e/d})$.
Pour $G'(d,d,2)$, si $e/2$ n'est pas une puissance d'un nombre premier,
alors $1 + \zeta^{-1}$ est inversible (car $-\zeta^{-1}$ est une racine de
l'unit\'e d'ordre soit $e$, soit $e/2$), et encore une fois il est clair
que $t = e/d$ pour $\bar v_1 = (1+\zeta^{-1})/(1-\zeta^{e/d})$. 

Enfin, si $e/2$ est bien une puissance d'un nombre premier, alors $d$ doit
\^etre une puissance du m\^eme nombre premier: \'ecrivons $e = 2p^a$ et $d
= p^b$, o\`u $b < a$ (puisqu'on a suppos\'e que $d < e/2$). Posons $c = a - b$;
alors on a que 
\[
\bar v_1 = \frac{1 + \zeta^{-1}}{1-\zeta^{2p^c}} = \frac{1+
\zeta^{-1}}{(1-\zeta^{p_c})( 1+ \zeta^{p_c})}. 
\]
Supposons d'abord que $p$ soit impair. Alors $-\zeta^{-1}$, $\zeta^{p^c}$,
$-\zeta^{p^c}$ sont des racines de l'unit\'e d'ordre $p^a$, $2p^b$, $p^b$,
respectivement. Le module de $\bar v_1$ (c'est-\`a-dire, le produit de
ses conjugu\'es par $\Gal(K/\Q)$) est donc 
\[
|\bar v_1| = \frac{|1+ \zeta^{-1}|}{|1-\zeta^{p^c}|\cdot |1+\zeta^{p^c}|} =
\frac{p}{1 \cdot p^{p^c}} = p^{1-p^c}. 
\]
(Voir, par exemple,~\cite{W}). On voit ici la n\'ecessit\'e d'avoir suppos\'e $c > 0$. Pour que $(1-\zeta^t)\bar v_1$ soit dans $\Z_K$, il faut que $|1-\zeta^t| \ge p^{p^c-1}$. Le plus petit tel $t$ est $t = 2p^c = e/d$, avec $|1-\zeta^t| =
p^{p^c}$. 

Dans le cas o\`u $p = 2$, le calcul est presque pareil: cette fois, les
ordres de $-\zeta^{-1}$, $\zeta^{p^c}$, $-\zeta^{p^c}$ sont $2^{a+1}$,
$2^{b+1}$, $2^{b+1}$, respectivement, et 
\[
|\bar v_1| = \frac{|1+ \zeta^{-1}|}{|1-\zeta^{p^c}|\cdot |1+\zeta^{p^c}|} = \frac{2}{2^{2^{c-1}} \cdot 2^{2^{c-1}}} = 2^{1-2^c}.
\]
Ensuite, le m\^eme argument montre que $t = e/d$ est le plus petit entier
tel que $(1-\zeta^t)\bar v_1 \in \Z_K$. 
\end{proof}

\end{document}